\documentclass[11pt, reqno]{amsart}
\usepackage{amsfonts,amssymb,latexsym,amsmath, amsxtra,url, mathrsfs}
\usepackage{verbatim}
\usepackage{enumitem}
\usepackage{tikz}
\usepackage{cancel}
\usepackage{bm}
\usepackage{bbm}

\theoremstyle{plain}
\newtheorem{theorem}{Theorem}[section]
\newtheorem*{theorem*}{Theorem}
\newtheorem{corollary}[theorem]{Corollary}

\newtheorem{lemma}[theorem]{Lemma}
\newtheorem{proposition}[theorem]{Proposition}
\newtheorem*{conjecture*}{Conjecture}

\theoremstyle{definition}

\newtheorem*{remark}{Remark}
\theoremstyle{remark}

\numberwithin{equation}{section}

\newcommand{\R}{\mathbb R}
\newcommand{\N}{\mathbb N}
\newcommand{\Z}{\mathbb Z}

\def\({\left(}
\def\){\right)}

\def\Log{\operatorname{Log}}

\def\lp{\left(}
\def\rp{\right)}

\def\l{\lambda}

\def\th{\theta}

\def\lp{\left(}
\def\rp{\right)}

\title{Statistics for unimodal sequences}
\author{Walter Bridges}
\author{Kathrin Bringmann}
\address{University of Cologne, Department of Mathematics and Computer Science, Weyertal 86-90, 50931 Cologne, Germany }
\email{wbridges@uni-koeln.de}
\email{kbringma@math.uni-koeln.de}

 

\begin{document}
\maketitle

\begin{abstract}
We prove a number of limiting distributions for statistics for unimodal sequences of positive integers by adapting a probabilistic framework for integer partitions introduced by Fristedt.  The difficulty in applying the direct analogue of Fristedt's techniques to unimodal sequences lies in the fact that the generating function for partitions is an infinite product, while that of unimodal sequences is not.  Essentially, we get around this by conditioning on the size of the largest part and working uniformly on contributing summands.  Our framework may be used to derive many distributions, and our results include joint distributions for largest parts and multiplicities of small parts.  We discuss ranks as well.  We further obtain analogous results for strongly unimodal sequences.
\end{abstract}

\section{Introduction and Statement of results}

A {\it partition} $\lambda$ of $n$ is a sequence of positive integers that sum to $n$,
$$
\lambda: \qquad \lambda_1 \geq \dots \geq \lambda_{\ell} > 0, \qquad \sum_{k=1}^{\ell} \lambda_k=n.
$$
We write $|\lambda|=n$ for the {\it size} of $\lambda$, set $p(n):= \#\{\lambda \vdash n\}$, and we define $p(0):=1$.  The generating function for partitions is the well-known infinite product
$$
P(q):=\sum_{\lambda} q^{|\lambda|}=\sum_{n \geq 0} p(n)q^n = \prod_{n \geq 1} \frac{1}{1-q^n}.
$$  An important result in the study of partition statistics is due to Erd\H{o}s and Lehner who proved that, as $n \to \infty$, the largest part of almost all partitions of $n$ is roughly $A\sqrt{n}\log(A\sqrt{n})$ and varies from this mean by an extreme value distribution.  Here and throughout the article, let $A:=\frac{\sqrt{6}}{\pi}$.
\begin{theorem}[Theorem 1.1 of \cite{EL}]\label{T:ErdosLehner}
For $v \in \mathbb{R}$, we have
$$
\lim_{n \to \infty} \frac{\#\left\{\lambda \vdash n : \frac{\lambda_1-A\sqrt{n}\log\left(A\sqrt{n}\right)}{A\sqrt{n}} \leq v \right\}}{p(n)} = e^{-e^{-v}}.
$$
\end{theorem}

Erd\H{o}s and Lehner's proof used only straightforward recurrences and the Hardy--Ramanujan asymptotic formula for $p(n)$ (\cite{HR}, equation (1.41)).  After a series of papers of Szalay--Tur\'{a}n and Erd\H{o}s--Tur\'{a}n (see \cite{ET,ST1,ST2,ST3}),  Fristedt introduced what has proved to be an indispensable probabilistic technique, allowing him to greatly extend previous distributions \cite{F}.  In particular, Theorem \ref{T:ErdosLehner} was extended to a joint distribution for the $t_n$ largest parts, where $t_n=o(n^{\frac 14})$.
\begin{theorem}[Theorems 2.5 and 2.6 of \cite{F}]\label{T:Fristedtlargeparts}
For any integer $t_n=o(n^{\frac 14})$ and $\{v_t\}_{t =1}^{t_n} \subset \mathbb{R}^{t_n}$, the following limit vanishes
$$
\lim_{n \to \infty} \left(\frac{\#\left\{\lambda \vdash n : \frac{\lambda_t-A\sqrt{n}\log\left(A\sqrt{n}\right)}{A\sqrt{n}}  \leq v_t \ \text{for $1 \leq t \leq t_n$} \right\}}{p(n)} 
- \int_{-\infty}^{v_1} \cdots \int_{-\infty}^{v_{t_n}} f(u_1,\dots, u_{t_n}) du_{t_n} \cdots du_1 \right),
$$
where 
$$
f(u_1,\dots, u_{t_n}):= \begin{cases} e^{-\sum_{t=1}^{t_n} u_t - e^{-u_{t_n}}} & \text{if $u_1 \geq \dots \geq u_{t_n}$,} \\ 0 & \text{otherwise.} \end{cases}
$$
\end{theorem}

\begin{remark}
The distribution on the right-hand side was interpreted in \cite{F} as a limiting Markov chain.  Another equivalent distribution was given by Pittel (\cite{Pi}, lemma on p. 127) in his proof of Wilf's conjecture that the limiting proportion of partitions of $n$ whose parts are the vertex degrees of a simple graph is 0.  In follow-up work, Pittel further improved on the above by proving explicit convergence rates \cite{P3}.  We leave a similar strengthening of our results below, as well as any analogous applications, as interesting open problems.
\end{remark}

\begin{remark}
Fristedt obtained a stronger version of the above which he stated in terms of the L\'{e}vy--Prokhorov distance between measures (see \cite{Bi1}, p. 72).  Our results below too could be strengthened in this way, but we prefer the simpler statements which show only convergence in distribution.
\end{remark}
Fristedt found a wide array of limiting distributions by introducing what is now known as a Boltzmann model.  This is now a standard technique for studying the statistical behavior and constructing sampling algorithms of many combinatorial structures (see \cite{DFLS}), but our exposition here is self-contained.  The Boltzmann model replaces the uniform probability measure on $\{\lambda \vdash n\}$ with a measure on all partitions $\lambda$, by defining
\begin{equation}\label{E:Boltzmannmodelpartitionsdef}
Q_q(\lambda):= \frac{q^{|\lambda|}}{P(q)}, \qquad \text{for $q \in (0,1)$.}
\end{equation}
Since $Q_q$ agrees on all $\lambda \vdash n$, the Boltzmann model, when conditioned on $\lambda \vdash n$, equals the uniform measure on partitions of $n$, and thus this technique is often called a {\it conditioning device}.  In statistical mechanics, the Boltzmann model is known as the macro-canonical ensemble, and the uniform measure on $\{\lambda \vdash n\}$ is known as the  micro-canonical ensemble (see \cite{V}).
Under $Q_q$, one directly gains independence of the relevant random variables that one lacks under the uniform probability measure.  This is precisely because $P(q)$ is an infinite product, and much of the work applying Boltzmann models to study statistics for partitions relies on product generating functions.  In this article, we show that a conditioned Boltzmann model may still be useful even without an infinite product generating function.\footnote{As pointed out to us by one of the referees there is a somewhat similar construction in work of Pittel \cite{P2} on partitions with bounded min/max ratio.}  We demonstrate this in the case of unimodal sequences.

A {\it unimodal sequence} $\lambda$ of size $n$ is a generalization of an integer partition, in which parts are allowed to increase and then decrease\footnote{For convenience, we reverse the usual convention of indexing so that here $\l_t^{[L]}, \l_{t+1}^{[L]} $ are decreasing in $t$.},
\begin{equation}\label{E:uniseqdef}
\lambda: \qquad \lambda_{r}^{[L]} \leq \dots \leq \lambda_{1}^{[L]} \leq \lambda_{\mathrm{PK}} \geq \lambda_{1}^{[R]} \dots \geq \lambda_{s}^{[R]} \qquad \text{and} \qquad \lambda_{\mathrm{PK}} + \sum_{k=1}^{r} \lambda_{k}^{[L]} + \sum_{k=1}^{s} \lambda_{k}^{[R]}  = n.
\end{equation}
We write $\mathcal{U}_n$ for the set of unimodal sequences of size $n$, with $\mathcal{U}_0$ containing a single empty sequence.  Let $u(n):= \#\mathcal{U}_n$ and $\mathcal{U}:=\cup_{n \geq 0} \mathcal{U}_n$.  The special part $\lambda_{\mathrm{PK}}$ is called the {\it peak}.  The generating function for $u(n)$ is obtained by summing over the size of peaks as
$$
U(q):=\sum_{\lambda \in \mathcal{U}} q^{|\lambda|} = \sum_{n \geq 0} u(n)q^n=1+\sum_{m \geq 1} q^m \prod_{k=1}^{m} \frac{1}{\left(1-q^k\right)^2}.
$$

	We note that there are two slightly different definitions of unimodal sequences in the literature which differ in whether or not the peak is specified.  The second author and Mahlburg studied the asymptotic behavior of various types of unimodal sequences, where what we call unimodal sequences in this paper were referred to as {\it stacks with summits} on p. 196 of \cite{BM}, as opposed to {\it stacks}, where the peak is unspecified.  (This is equivalent to forcing the inequality $\lambda_{\mathrm{PK}}\geq\lambda_{s}^{[R]}$ to be strict.)  All of our results here hold for these stacks as well; it is simply easier to state and prove our results with the above definition. Stacks seem to have been introduced by Auluck in \cite{Au}, who called them {\it type B partitions}, whereas stacks with summits are also called {\it V-partitions} in $\S 2.5$ of \cite{St}.  Andrews \cite{A1} has also called these sequences {\it convex compositions} when the peak is strictly larger than the other parts (which is equivalent to adding 1 to the peak of unimodal sequences). As pointed out to us by one of the referees statistics for unimodal sequences were also studied in the context of statistical mechanics by Temperley (see \cite{T}, Section 3), who gave a heuristic for the limiting shape of diagrams.  Limit shapes for several types of unimodal sequences were rigorously proved by the first author in \cite{B21}.  Recent work of the second author, Jennings-Shaffer, Mahlburg and Rhoades \cite{BJSM,BJSMR} proved the distribution of the rank, defined below, a statistic connected to modular forms.

To state our results, let $\mathbf{P}_{n}$ be the uniform probability measure on $\mathcal{U}_n$, and introduce the following random variables on $\mathcal{U}$:
\begin{itemize}[leftmargin=*]
\item  Let ${\rm PK}(\lambda):=\lambda_{\mathrm{PK}}$ denote the peak of a sequence.
\item  Let $X_{k}^{[L]}(\lambda)$ (resp. $X_{k}^{[R]}(\lambda)$) denote the number of parts in $\lambda$ equal to $k$ and to the left (resp. right) of the peak.
\item  Let $Y_{t}^{[L]}(\lambda):=\lambda_t^{[L]}$ (resp. $Y_{t}^{[R]}$) denote the $t$-th largest part in $\lambda$ to the left (resp. right) of the peak.
\item  Let $N(\lambda):=\sum_{k \geq 1} k(X_{k}^{[L]}(\lambda) + X_{k}^{[R]}(\lambda)) + {\rm PK}(\lambda)$ denote the {\it size} of $\lambda$. 
\end{itemize}

\noindent{\bf Example.}  If $\lambda$ is the unimodal sequence $1+1+2+3+3+{\bf 3}+1$ (the peak is boldfaced), then 
\begin{align*}
&X_{1}^{[L]}(\lambda)=2, \ X_{2}^{[L]}(\lambda)=1, \ X_{3}^{[L]}(\lambda)=2, \ {\rm PK}(\lambda)=3, \  X_{3}^{[R]}(\lambda)=0,\  X_{2}^{[R]}(\lambda)=0, \ X_{1}^{[R]}=1,\\  
&N(\lambda)=14.
\end{align*}

Note that the $Y_t^{[L]}$ are explicitly determined in terms of the $X_k^{[L]}$ as
$$
Y_t^{[L]}:=\sup\left\{\ell : \sum_{k \geq \ell} X_k^{[L]} \geq t \right\}.
$$
A similar statement holds for $Y_t^{[R]}$.

Our first result is an analogue of Theorem \ref{T:ErdosLehner}.  Here and throughout the article, let $B:=\frac{\sqrt{3}}{\pi}$.
\begin{theorem}\label{T:PKdist}
For $v \in \mathbb{R}$, we have
$$
\lim_{n \to \infty} \mathbf{P}_{n}\left(\frac{{\rm PK}-B\sqrt{n}\log\left(2B\sqrt{n}\right)}{B\sqrt{n}} \leq v \right)=e^{-e^{-v}}.
$$
Moreover, if $\mathbf{E}_n$ denotes expectation under $\mathbf{P}_{n}$, then we have
$$
\mathbf{E}_n({\rm PK}) = B\sqrt{n}\log\left(2B\sqrt{n}\right)+B\gamma \sqrt{n}(1+o(1)),
$$
where $\gamma$ is the Euler--Mascheroni constant.
\end{theorem}

\begin{remark} Note that $\mathcal{U}_n$ is almost the same as the set of pairs of partitions of total size $n$, so a reasonable heuristic guess would be that  Theorem \ref{T:ErdosLehner} holds with $\lambda_1 \mapsto {\rm PK}$ and $n \mapsto \frac{n}{2}$.  However, the presence of the extra factor of $\log(2)$ is justified by Theorem \ref{T:largeparts} below.
\end{remark}

\begin{remark}
It would be interesting to improve the error term for the mean, as Ngo and Rhoades did for Theorem 1.1 \cite{NR}.  We discuss this further in Section \ref{S:moments}.
\end{remark}

We also prove the following analogue of Theorem \ref{T:Fristedtlargeparts}.
\begin{theorem}\label{T:largeparts}
For any integer $t_n=o(n^{\frac 14})$ and $\{v_t^{[j]}\}_{1 \leq t \leq t_n, j \in \{L,R\}} \subset \mathbb{R}^{2t_n}$, the following difference vanishes as $n \to \infty$,
\begin{multline*} 
	\mathbf{P}_{n}\left(\frac{{\rm PK}-B\sqrt{n}\log\left(2B\sqrt{n}\right)}{B\sqrt{n}} \leq v_0, \ \frac{Y_{t}^{[j]}-B\sqrt{n}\log\left(2B\sqrt{n}\right)}{B\sqrt{n}} \leq v_{t}^{[j]}, \ \ \text{for $j \in \{L,R\}$, \ $1 \leq t \leq t_n$}\right) \\
	- \int_{-\infty}^{v_0} \int_{-\infty}^{v_{1}^{[L]}} \int_{-\infty}^{v_{1}^{[R]}} \cdots \int_{-\infty}^{v_{t_n}^{[L]}} \int_{-\infty}^{v_{t_n}^{[R]}} F\left(u_0,u_{1}^{[L]},u_{1}^{[R]}, \dots, u_{t_n}^{[L]},u_{t_n}^{[R]}\right) du_{t_n}^{[R]} \cdots du_0,
\end{multline*}
where
\begin{multline*}
F\left(u_0,u_{1}^{[L]},u_{1}^{[R]}, \dots, u_{t_n}^{[L]},u_{t_n}^{[R]}\right) \\:= \begin{cases}\frac{1}{2^{2t_n}}e^{-u_0-\sum_{1 \leq t \leq t_n} \left(u_t^{[L]}+u_t^{[R]}\right) - \frac{e^{-u_{t_n}^{[L]}}}{2} - \frac{e^{-u_{t_n}^{[R]}}}{2}} & \text{if $u_0 \geq u_1^{[j]} \geq \dots \geq u_{t_n}^{[j]}$ for $j \in \{L,R\}$}, \\ 0  & \text{otherwise.}   \end{cases}
\end{multline*}
\end{theorem}

Next, we show that the joint distribution of the numbers of small parts, when re-scaled,  behaves as the joint distribution for independent and identically distributed exponential random variables, and that of medium sized parts behaves geometrically; this is an analogue of Theorem 2.2 of \cite{F}.
\begin{theorem}\label{T:smallparts}
Let $k_n=o(n^{\frac 14})$ be an integer and $\{v_t^{[j]}\}_{1 \leq k \leq k_n, j \in \{L,R\}} \subset [0,\infty)^{2k_n}$. Then
$$
\lim_{n \to \infty} \left(\mathbf{P}_{n}\left(\frac{kX_{k}^{[j]}}{B\sqrt{n}} \leq v_{k}^{[j]}, \ \ \text{for $1 \leq k \leq k_n$ and $j \in \{L,R\}$} \right) - \prod_{\substack{1 \leq k \leq k_n \\ j \in \{L,R\}}} \int_{0}^{v_{k}^{[j]}} e^{-u_{k}^{[j]}}du_{k}^{[j]} \right)=0.
$$
Moreover, for an integer $k=o(n^{\frac{1}{2}})$ and $v^{[L]}, v^{[R]} \in \mathbb{R}$, we have
$$
\lim_{n \to \infty} \mathbf{P}_n\left(\frac{kX_k^{[L]}}{B\sqrt{n}} \leq v^{[L]}, \frac{kX_k^{[R]}}{B\sqrt{n}} \leq v^{[R]} \right) = \left(1-e^{-v^{[L]}}\right)\left(1-e^{-v^{[R]}}\right),
$$
and for $k=\lfloor cn^{\frac{1}{2}} \rfloor$ with $v^{[L]},v^{[R]} \in \mathbb{N}_0$, we have
$$
\lim_{n \to \infty} \mathbf{P}_n\left(X_k^{[L]} \leq v^{[L]}, X_k^{[R]} \leq v^{[R]} \right) = \left(1-e^{-\frac{c}{B}\left(v^{[L]}+1\right)}\right)\left(1-e^{-\frac{c}{B}\left(v^{[R]}+1\right)}\right).
$$
\end{theorem}

It is combinatorially obvious that $\mathbf{P}_{n}(X_{k}^{[L]} \leq X_{k}^{[R]})=\frac{1}{2}$, simply by swapping left and right parts on every element of $\mathcal{U}_n$.  A more refined measure of symmetry in the small parts would be a joint distribution for the differences of $X_{k}^{[L]}$ and $X_{k}^{[R]}$ for $k \leq k_n$.  As a corollary of the above, we show that this behaves like the joint distribution of independent Laplace distributions.  Below and throughout, $\chi_S:=1$ if a statement $S$ holds and 0 otherwise. 
\begin{corollary}\label{C:smallpartsskewness}
Let $k_n=o(n^{\frac 14})$ be an integer and let $\{v_k\}_{1 \leq k \leq k_n} \subset \mathbb{R}^{k_n}$.  Then the following limit vanishes, 
\begin{align*}
&\lim_{n \to \infty} \left(\mathbf{P}_{n}\left(\frac{k\left(X_{k}^{[L]}-X_{k}^{[R]}\right)}{B\sqrt{n}} \leq v_{k}, \ \ \text{for $1 \leq k \leq k_n$} \right) - \prod_{1 \leq k \leq k_n} \left(\left(1-\frac{e^{-v_k}}{2}\right)\chi_{v_k \geq 0} + \frac{e^{v_k}}{2} \chi_{v_k < 0} \right)  \right).
\end{align*}
\end{corollary}
By summing the differences $X_{k}^{[L]}-X_{k}^{[R]}$ over all $k$, we obtain the {\it rank} of a unimodal sequence,
\begin{equation}\label{E:rankdef}
{\rm rank}(\lambda):=\sum_{k \geq 1}\left( X_{k}^{[L]}(\lambda)-X_{k}^{[R]}(\lambda) \right).
\end{equation}
Bringmann, Jennings-Shaffer, and Mahlburg proved that at the scaling of $\sqrt{n}$ the rank obeys a logistic distribution.
\begin{theorem}[Proposition 1.2 of \cite{BJSM}]\label{T:rank}
For $x \in \mathbb{R}$, we have
$$
\lim_{n \to \infty} \mathbf{P}_{n}\left(\frac{{\rm rank}}{B\sqrt{n}} \leq x \right) = \frac{1}{1+e^{-x}}.
$$
\end{theorem}
Bringmann, Jennings-Shaffer, and Mahlburg used the method of moments, an approach that relies on suitable two-variable generating functions and is independent of the present work.  We prove the following related result.
\begin{theorem}\label{T:totalsmallparts}
For any integer $k_n=o (n^{\frac{1}{2}})$ with $k_n \to \infty$ and $v_L,v_R \in \mathbb{R}$,
$$\lim_{n \to \infty} \mathbf{P}_{n}\left( \sum_{k \leq k_n} \frac{X_{k}^{[j]}-B\sqrt{n}\log(k_n)}{B\sqrt{n}} \leq v_{j}, \ \text{for $j \in \{L,R\}$}\right) = e^{-e^{-v_{L}}-e^{-v_R}}.$$
\end{theorem}
Thus, the total small part counts on the left and the right behave as independent extreme value distributions when the mean of $B\sqrt{n}\log (k_n)$ is subtracted, and their convolution gives the logistic distribution in Theorem \ref{T:rank}.  Our techniques are not robust enough to prove Theorem \ref{T:rank}, but when combined with the following, Theorem \ref{T:totalsmallparts} is highly suggestive of Theorem \ref{T:rank}.
\begin{proposition}[see Theorem 1.2 of \cite{B21}]\label{P:totallargeparts}
For any fixed $\varepsilon > 0$ and any fixed $a >0$, we have
$$
\lim_{n \to \infty} \mathbf{P}_{n}\left(\left|\frac{1}{\sqrt{n}} \sum_{k \geq a\sqrt{n}} \left(X_{k}^{[L]}-X_{k}^{[R]}\right)\right| < \varepsilon \right) = 1.
$$
\end{proposition}
Note that our methods also easily extend to the so-called $t$-th successive ranks, defined by limiting the range in \eqref{E:rankdef} to $k \geq t$, whereas extension of the techniques in \cite{BJSM} appears to be non-trivial.

Finally, we can prove results similar to Theorems \ref{T:largeparts} and \ref{T:smallparts} for strongly unimodal sequences, defined so that the inequalities in \eqref{E:uniseqdef} are strict.  Let $\mathcal{U}_n^*$ denote the set of strongly unimodal sequences of size $n$, let $u^*(n):=\#\mathcal{U}_n^*$ and let $\mathbf{P}_{n}^*$ denote the uniform probability measure on $\mathcal{U}^*_n$.

\begin{theorem}\label{T:*largesmallparts}
For strongly unimodal sequences, Theorems \ref{T:PKdist} and \ref{T:largeparts} hold with $\mathbf{P}_{n}$ replaced by $\mathbf{P}_{n}^*$ and $B= \frac{\sqrt{3}}{\pi}$ replaced by $A=\frac{\sqrt{6}}{\pi}$.
\end{theorem}
Since a part can occur at most once on either side of the peak in a strongly unimodal sequence, it is natural to expect the following analogue of Theorem 9.2 of \cite{F}.  Note that here we can take a larger range for $k_n$ than in Theorem \ref{T:smallparts}.  Here and throughout, boldface letters represent vectors (except when we use them for probability measures).
\begin{theorem}\label{T:*smallparts}
Suppose that $k_n=o(n^{\frac{1}{2}})$ and let ${\bm w} \in \{0,1\}^{2k_n}$.  Then
$$
\lim_{n \to \infty} \left(\mathbf{P}_{n}^*\left( \left(X_{k}^{[j]}\right)_{j \in \{L,R\}, 1 \leq k \leq k_n} = {\bm w} \right) - \frac{1}{2^{2k_n}} \right)=0.
$$
\end{theorem}

We content ourselves  with the above results, as once our machinery is described many of the details follow \cite{F} closely.  It would be interesting to apply our methods to the other types of stacks and unimodal sequences discussed in \cite{BM} or to extend the work of Pittel \cite{P} to unimodal sequences, where Fristedt's methods were used to study distributions of mid-range parts in partitions.

Section \ref{S:Preliminaries} contains some preliminaries needed for asymptotic analysis. In Section \ref{S:Fristedt}, we briefly recall Fristedt's machinery and in Section \ref{S:Boltzmann}, we introduce our modification, a conditioned Boltzmann model for unimodal sequences.  We believe these methods should be useful in contexts beyond the present work, especially for non-product generating functions in the theory of partitions.  Section \ref{S:details} contains the details of the proofs of our main theorems.  Lastly, in Section \ref{S:moments} we discuss the moment generating functions for ${\rm PK}$, a method independent from Boltzmann models, but one which may lead to more precise results for ${\rm PK}$.

\section*{Acknowledgements} The first author is partially supported by the supported by the SFB/TRR 191 ``Symplectic Structures in Geometry, Algebra and Dynamics'', funded by the DFG (Projektnummer 281071066 TRR 191).  The second author has received funding from the European Research Council (ERC) under the European Union's Horizon 2020 research and innovation programme (grant agreement No. 101001179). We thank the referees for their useful comments which improved the exposition.

\section{Preliminaries}\label{S:Preliminaries}

\subsection{Notation}

We use the following standard asymptotic notation.  We write
\begin{itemize}[leftmargin=*]
\item $f(n) \sim g(n)$ if $\lim_{n \to \infty} \frac{f(n)}{g(n)} =1$,
\item  $f(n)=O(g(n))$ (or $f(n) \ll g(n)$) if the quotient $\frac{f(n)}{g(n)}$ is bounded as $n \to \infty$,
\item $f(n)=o(g(n))$ if $\lim_{n \to \infty} \frac{f(n)}{g(n)}=0$,
\item  $f(n)=\omega(g(n))$ if the quotient $\frac{f(n)}{g(n)}$ is unbounded as $n \to \infty$, and
\item  $f(n) \asymp g(n)$ if $f(n)=O(g(n))$ and $g(n)=O(f(n))$.
\end{itemize}

To concisely write products in generating functions, we employ the standard $q$-factorial notation, defined for $n \in \N_0 \cup \{\infty\}$ by
$$
(a)_n=(a;q)_n:=\prod_{j=0}^{n-1} \left(1-aq^j\right).
$$

\subsection{Euler--Maclaurin summation, logarithmic series, and integral calculations}

We need the following two variants of Euler--Maclaurin summation. Let $\{x\}:=x-\lfloor x \rfloor$ denote the fractional part of $x$. Note that the first variant is recovered by taking $a \to 1^-$ in Theorem B.5 of \cite{MV}.

\begin{lemma}[Theorem B.5 in \cite{MV}] \label{L:eulermac}
For $N \in \mathbb{N}$ and continuously differentiable $g:\mathbb{R} \to \mathbb{C}$, we have 
\begin{align}
\sum_{k=1}^N g(k) &= \int_1^N g(u)du + \frac{1}{2}\left(g(N)+g(1)\right) + \int_1^N  \left( \{u\} - \frac{1}{2}\right) g'(u)du \label{E:eulermac1} \\
&= \int_0^N g(u)du + \frac{1}{2}\left(g(N)-g(0)\right) + \int_0^N  \left( \{u\} - \frac{1}{2}\right) g'(u) du. \label{E:eulermac2} 
\end{align}
\end{lemma}

The following lemma is useful for the approximation logarithmic series by Taylor expansions.  Here and throughout, we write $\Log$ for the principal branch of the complex logarithm, and  $\log$ when $\Log$ is restricted to the positive real axis.
\begin{lemma}[Lemma 1 in \cite{R} with Lemma 1 in \cite{B20}]\label{L:Romik}
There exists a constant $C>0$ such that for all $0<x<1$ and $s \in \R$, we have
$$\left|\Log\left(\frac{1\pm x}{1\pm xe^{is}}\right)-\frac{isx}{1\pm x}+ \frac{s^2x}{2(1\pm x)^2} \right| \leq C \frac{x|s|^3}{(1-x)^3}.$$
\end{lemma}

We also need the following lemma concerning the asymptotic behavior of a certain product.  This should be compared with a similar formula on p. 723 of \cite{F}.

\begin{lemma}\label{L:largepartsproduct}
Uniformly in $v \geq - \frac{\log (n)}{8}$ as $n \to \infty$, we have
\begin{align}
 \prod_{k > A\sqrt{n}\left(v+\log(A\sqrt{n})\right)} \left(1-e^{-\frac{k}{A\sqrt{n}}}\right) &\sim e^{-e^{-v}},\label{E:largepartsproduct1}\\
 \prod_{k > A\sqrt{n}\left(v+\log(A\sqrt{n})\right)} \frac{1}{1+e^{-\frac{k}{A\sqrt{n}}}} &\sim e^{-e^{-v}}.\label{E:largepartsproduct2}
\end{align}
\end{lemma}
\begin{proof}
Equation (\ref{E:largepartsproduct1}) is proved in equation (6.10) of \cite{F}. Next note that \eqref{E:largepartsproduct2} is equivalent to
$$
-\lim_{n \to \infty} \sum_{k > A\sqrt{n}\left(v+\log(A\sqrt{n})\right)} \log\left(1+e^{-\frac{k}{A\sqrt{n}}}\right) +e^{-v}= 0,
$$
uniformly in $v \geq -\frac{\log(n)}{8}$.  After using the inequality $-x < - \log(1+x)<-x+\frac{x^2}{2}$ for $x \in (0,1)$, a short calculation shows
$$
O\left(e^{-v-\frac{\log (n)}{2}}\right) <- \sum_{k > A\sqrt{n}\left(v+\log(A\sqrt{n})\right)} \log\left(1+e^{-\frac{k}{A\sqrt{n}}}\right) +e^{-v} <O\left(e^{-v-\frac{\log(n)}{2}}\right) + O\left(e^{-2v-\frac{\log(n)}{2}}\right).
$$  Both terms are $o(1)$ if $v\geq-\frac{\log(n)}{8}$, completing the proof.
\end{proof}

Another integral evaluation gives the distribution of the sum of independent exponentially distributed random variables with means $1$, $\frac{1}{2}$, $\dots$, $\frac{1}{n}$; the proof follows using induction.
\begin{lemma}\label{L:stochYuleintegral}
For $x \geq 0$ and for $n \in \mathbb{N}$, we have
$$
n!\int_0^x \int_0^{x-u_n} \cdots \int_0^{x-u_n-\dots - u_2} e^{-\sum_{i=1}^{n} ju_j} du_1 \cdots du_n = \left(1-e^{-x}\right)^n.
$$
\end{lemma}

\subsection{Saddle-point method}

We require the following specific case of the saddle-point method for evaluating Cauchy integrals.  Our variant is essentially the one used by Fristedt \cite{F} in the proof of Proposition 4.5, and also in the proofs of \cite{R} Proposition 3 and \cite{B20} Proposition 3, although these are all stated in probabilistic terminology.  This may also be compared with \cite{FS} Chapter VIII, Figure VIII.4.
\begin{proposition}\label{P:saddlepointmethod}
Let $\{g_n\}_{n \geq 1}$ be a sequence of twice continuously differentiable functions.  Suppose that for all sufficiently small fixed $\varepsilon>0$, after decomposing the integral as
$$
\int_{-\frac{1}{2}}^{\frac{1}{2}} \exp\left(g_n(2\pi i \theta)\right)d \theta = \int_{-\varepsilon n^{-\frac{1}{2}}}^{\varepsilon n^{-\frac{1}{2}}} \exp\left(g_n(2\pi i \theta)\right)d \theta + \int_{\varepsilon n^{-\frac{1}{2}}<|\theta| \leq \frac{1}{2}} \exp\left(g_n(2\pi i \theta)\right)d \theta,
$$
the following hold as $n \to \infty$.
\begin{enumerate}[leftmargin=*,label=\rm{(\arabic*)}]
\item We have that $g_n(0)\asymp n^{\frac{1}{2}}$ and $g_n''(0) \asymp n^{\frac{3}{2}}$ where the implied constants are positive real numbers, $g_n'(0) = o(n^{\frac{3}{4}})$, and we have a quadratic approximation on the ``major arc'': for a constant independent of $\varepsilon$, we have 
$$
\left|g_n(2\pi i \theta)-g_n(0)-g_n'(0)2\pi i \theta - g_n''(0) \frac{(2\pi i \theta)^2}{2} \right| = O\left(\theta^3n^2\right), \qquad \text{for $|\theta|\leq \varepsilon n^{-\frac{1}{2}}$.}
$$
\item The ``minor arc'' is negligible: 
$$\lim\sup_{n \to \infty} \frac{\mathrm{Re}\left(g_n(2\pi i \theta)\right)-g_n(0)}{n^{\frac{1}{2}}} < -\delta_{\varepsilon}, \qquad \text{for some $\delta_{\varepsilon} >0$ and $\varepsilon n^{-\frac{1}{2}} |\theta|\leq \frac{1}{2}$}.$$
\end{enumerate}
Then
$$
\int_{-\frac{1}{2}}^{\frac{1}{2}} \exp\left(g_n(2 \pi i \theta)\right)d\theta \sim \frac{e^{g_n(0)}}{\sqrt{2\pi g_n''(0)}}.
$$
\end{proposition}

\begin{proof}
If we prove that 
$$
\int_{-\varepsilon n^{-\frac{1}{2}}}^{\varepsilon n^{-\frac{1}{2}}} \exp\left(g_n(2 \pi i \theta)\right)d\theta \sim \frac{e^{g_n(0)}}{\sqrt{2\pi g_n''(0)}},
$$
then the proposition follows by assumption (2).  Thus we focus on the major arc and use assumption (1) to write
$$
\int_{-\varepsilon n^{-\frac{1}{2}}}^{\varepsilon n^{-\frac{1}{2}}} \exp\left(g_n(2 \pi i \theta)\right)d\theta
= e^{g_n(0)}\int_{-\varepsilon n^{-\frac{1}{2}}}^{\varepsilon n^{-\frac{1}{2}}} \exp\left( g_n'(0) 2\pi i \theta + g_n''(0) \frac{(2 \pi i \theta)^2}{2}  + O\left(n^2\theta^3\right)\right)d \theta.
$$
We then substitute $\theta \mapsto \frac{\theta}{\pi \sqrt{2g_n''(0)}}$ to get
\begin{align*}
&\frac{e^{g_n(0)}}{\pi \sqrt{2g_n''(0)}}\int_{-\pi\varepsilon \sqrt{2g_n''(0)} n^{-\frac{1}{2}}}^{\pi\varepsilon \sqrt{2g_n''(0)} n^{-\frac{1}{2}}} \exp\left(\frac{g_n'(0)}{\pi \sqrt{2g''(0)}} 2 \pi i \theta - \theta^2\left(1+O\left(\frac{\theta}{g''_n(0)^{\frac{3}{2}}} n^2\right)\right)\right)d\theta \\
&\sim \frac{e^{g_n(0)}}{\pi \sqrt{2g_n''(0)}} \int_{-\pi\varepsilon \sqrt{2g_n''(0)} n^{-\frac{1}{2}}}^{\pi\varepsilon \sqrt{2g_n''(0)} n^{-\frac{1}{2}}} \exp\left(-\theta^2\left(1+O\left(\theta n^{-\frac{1}{4}}\right)\right)\right)d\theta,
\end{align*}
where by hypothesis the constant is independent of $\varepsilon$.  But $\pi \sqrt{2g_n''(0)} n^{-\frac{1}{2}} \asymp n^{\frac{1}{4}}$, and hence $\varepsilon$ can be taken small enough so that $1+O(\theta n^{\frac{1}{4}})$ is bounded by a positive constant.  Since the integrand converges pointwise in $\theta$ to $e^{-\theta^2}$, Lebesgue's Dominated Convergence Theorem gives
$$
\frac{e^{g_n(0)}}{\pi \sqrt{2g_n''(0)}} \int_{-\pi\varepsilon \sqrt{2g_n''(0)} n^{-\frac{1}{2}}}^{\pi\varepsilon \sqrt{2g_n''(0)} n^{-\frac{1}{2}}} \exp\left(-\theta^2\left(1+O\left(\theta n^{-\frac{1}{4}}\right)\right)\right)d\theta\sim \frac{\sqrt{\pi}e^{g_n(0)}}{\pi \sqrt{2g_n''(0)}}=\frac{e^{g_n(0)}}{\sqrt{2\pi g_n''(0)}},
$$
as claimed.
\end{proof}

\subsection{Probability theory}

We require only the notions of random variables and their distributions, as well as the total variation metric, $d_{\mathrm{TV}}$, which is defined on measures $\mu$ and $\nu$ on $\mathbb{R}^{d}$ by
$$
d_{\mathrm{TV}}(\mu,\nu):=\sup_{\substack{B \subset \mathbb{R}^{d} \\ \text{Borel}}} \left(\mu(B)-\nu(B)\right).
$$
The Borel sets are the $\sigma$-algebra on $\mathbb{R}^{d}$ generated by open sets.  We also use Chebyshev's Inequality. 
\begin{lemma}[Chebyshev's Inequality]\label{L:Chebyshev}
If $X$ is a square-integrable random variable under a probability measure $\mathbf{P}$ with finite expectation $m$ and variance $\sigma$, then for any $t>0$
$$
\mathbf{P}(|X-m| \geq t) \leq \frac{\sigma^2}{t^2}.
$$
\end{lemma}

\section{Boltzmann models for partitions and the work of Fristedt}\label{S:Fristedt}

Here we recall the primary methods of Fristedt in Section 4 of \cite{F}. Let $P_n$ be the uniform probability measure on partitions of size $n$.  In analogy to our situation, let $X_k$ be the random variable giving the number of $k$'s in a partition and let $N=\sum_{k \geq 1} kX_k$ denote the size of a partition.  The key point is that the $X_k$ are independent under $Q_q$ (but not under $P_n$).

\begin{proposition}[Proposition 4.1 of \cite{F}]\label{P:Fristedt4.1}
	The $X_k$ are independent under $Q_q$.  Furthermore, $Q_q(X_k=\ell)=q^{k\ell}(1-q^k)$, so $X_k$ is geometric with mean $\frac{q^k}{1-q^k}$.
\end{proposition}

By definition $Q_q$ agrees on all $\lambda \vdash n$, so an immediate consequence is that,
$$
P_n=Q_q( \cdot | N=n),
$$
i.e., $Q_q$, when conditioned on $N=n$, agrees with $P_n$.  (This can also be shown directly from the definition of $Q_q$.)  Now if $W_n:\mathcal{P} \to \mathbb{R}^{d_n}$ is a random vector on partitions which is defined in terms of the $X_k$, the probability measure $Q_q$ is used to bridge the gap between the probability distribtion $P_n(W_n^{-1})$ on $\mathbb{R}^{d_n}$ and some explicitly given distributions $\nu_n$ on $\mathbb{R}^{d_n}$.  For a particular $q=q(n)$, Fristedt showed that $P_n(W_n^{-1})$ and $Q_q(W_n^{-1})$ are close in the sense of total variation given a very general condition on $W_n$.

\begin{proposition}[Proposition 4.6 of \cite{F}]\label{P:Fristedt4.6}
Let $K_n$ be a set of integers such that $W_n:\mathcal{P} \to \mathbb{R}^{d_n}$ is determined by $X_k$ for $k \in K_n$ with probability 1.  Let $q=q(n)=e^{-\frac{\pi}{\sqrt{6n}}}$.  If
$$
\sum_{k \in K_n} \frac{k^2q^k}{\left(1-q^k\right)^2}=o\left(n^{\frac 32}\right),
$$
then we have $d_{\mathrm{TV}}(P_n(W_n^{-1}), Q_q(W_n^{-1})) \to 0$.
\end{proposition}

The individual probabilities $Q_q(W_n={\bm w})$ are often straightforward to calculate, and Fristedt found limiting distributions for $W_n$ under $Q_q$ for a number of different $W_n$.  Together with Proposition \ref{P:Fristedt4.6}, this proved his many results.  The generality of Proposition \ref{P:Fristedt4.6} is the primary advantage of the probabilistic approach over a classical, direct Circle Method/saddle-point method approach.

In the following section, we prove an analogue of Proposition \ref{P:Fristedt4.6} for a conditioned Boltzmann model for unimodal sequences.  Our exposition is self-contained, but the interested reader is invited to examine Section 4 of \cite{F} for more probabilistic intuition for why Proposition \ref{P:Fristedt4.6} is true.

\section{Conditioned Boltzmann model}\label{S:Boltzmann}

To bridge the gap between $\mathbf{P}_{n}$ and the distributions in our main theorems, we introduce the Boltzmann model.  This probability measure is defined for any $q \in (0,1)$ on all of $\mathcal{U}$ by 
$$
\mathbf{Q}_{q}(\lambda):= \frac{q^{|\lambda|}}{U(q)},
$$
where $U(q):=1+\sum_{m \geq 1} \frac{q^m}{(q)_m^2}$ (see \cite{BM}, equation (1.8)).  As written, $\mathbf{Q}_q$ is not very useful for us because there is not a simple expression for the individual probabilities, $\mathbf{Q}_{q}(X_{k}^{[L]}=\ell)$, and in fact the $X_{k}^{[L]}$ and $X_{\ell}^{[R]}$ are not independent.  This lies in the fact that $U(q)$ is not a product.  However, the $m$-th summand in $U(q)$ is of course the product
$$
q^m\prod_{k=1}^{m}\frac{1}{\left(1-q^k\right)^2}=\sum_{\substack{\lambda \\ {\rm PK}(\lambda)=m}} q^{|\lambda|}.
$$
By conditioning $\mathbf{Q}_q$ on the event ${\rm PK}=m$, we do gain tractable expressions for the individual probabilities of $X_{k}^{[L]}=\ell$, and we can use the full power of Fristedt's techniques in \cite{F}.  Furthermore - and most importantly - we can do this uniformly for $m$ in the contributing range, and thus we are able to piece together the local distributions to obtain our global results.  We conclude this section with Proposition \ref{P:W_ncondition}, a direct analogue of Proposition \ref{P:Fristedt4.6}.

We do the analogous thing for strongly unimodal sequences, setting
$$
\mathbf{Q}^*_q(\lambda):=\frac{q^{|\lambda|}}{U^*(q)},
$$
where $U^*(q):=\sum_{m \geq 0} (-q)^2_{m-1}q^m$ is the generating function for strongly unimodal sequences (see \cite{A1}, where $u^*(n)=x_d(n)$). Throughout the article, we place $*$ in the superscript (and sometimes in the subscript) for the corresponding sets and functions defined on strongly unimodal sequences, calling attention to any significant differences when they arise.

Set $\mathbf{Q}_{q,m}:=\mathbf{Q}_q\left( \cdot | {\rm PK}=m \right)$ and $\mathbf{P}_{n,m}:=\mathbf{P}_{n} \left( \cdot | {\rm PK}=m \right)$.  Let $\mathcal{U}_{n,m} \subset \mathcal{U}_n$ be those sequences with peak $m$, and set $u_m(n):=\# \mathcal{U}_{n,m}$. For strongly unimodal sequences, it is more convenient for indexing to set $\mathbf{P}_{n,m}^*:=\mathbf{P}_{n}^* \left( \cdot | {\rm PK}=m+1 \right)$ and $\mathbf{Q}_{q,m}^*:=\mathbf{Q}_q^*\left( \cdot | {\rm PK}=m+1\right)$.

The following two lemmas are easily verified directly from the definitions.
\begin{lemma}\label{L:Q_{q,m}independenceetc}
\hspace{1mm}\newline
{\rm (1)} We have 
$$\mathbf{Q}_{q,m}(\lambda)=\begin{cases} (q)_m^2 q^{|\lambda|-m} & \text{if ${\rm PK}(\lambda)=m$,} \\ 0 & \text{otherwise.} \end{cases}$$
{\rm (2)} The set $\{X_{k}^{[j]} \}_{\substack{j \in \{L,K\} \\ k \geq 1}}$ is a set of independent random variables under $\mathbf{Q}_{q,m}$ with probability 

densities
$$
\mathbf{Q}_{q,m}\left(X_{k}^{[j]}=\ell\right)=\begin{cases} \left(1-q^k\right)q^{\ell k} & \text{if $k \leq m$,}  \\ 0 & \text{otherwise.} \end{cases}
$$

In particular,
$$
{\bf Q}_{q,m}(N=n)= \left[\zeta^{n-m}\right] \frac{(q;q)_m^2}{(\zeta q; \zeta q)_m^2}.
$$
{\rm (3)} We have $\mathbf{P}_{n,m}=\mathbf{Q}_{q,m}\left( \cdot | N=n \right)$.
\end{lemma}

For strongly unimodal sequences, we have the following analogue.

\begin{lemma}\label{L:*Q_{q,m}independenceetc} 
\hspace{1mm}\newline
{\rm ($1$)} We have
$$\mathbf{Q}_{q,m}^*(\lambda)=\begin{cases} (-q)_m^2 q^{|\lambda|-m-1} & \text{if ${\rm PK}(\lambda)=m+1$,} \\ 0 & \text{otherwise.} \end{cases}$$
{\rm ($2$)} The set $\{X_{k}^{[j]} \}_{\substack{j \in \{L,R\} \\ k \geq 1}}$ is a set of independent random variables under $\mathbf{Q}_{q,m}$ with probability 

densities
$$
\mathbf{Q}_{q,m}^*\left(X_{k}^{[j]}=\ell\right)=\begin{cases} \frac{1}{1+q^k} & \text{if $k \leq m$ and $\ell=0$},  \\ \frac{q^k}{1+q^k} & \text{if $k \leq m$ and $\ell=1$},  \\ 0 & \text{otherwise.} \end{cases}
$$

In particular,
$$
 {\bf Q}^*_{q,m}(N=n)= \left[\zeta^{n-m}\right] \frac{(-\zeta q; \zeta q)_m^2}{(-q;q)_m^2}
$$
{\rm ($3$)} We have $\mathbf{P}_{n,m}^*=\mathbf{Q}_{q,m}^*\left( \cdot | N=n \right)$.
\end{lemma}

Heuristically, it makes sense to choose $q=q(n)\in(0,1)$ to maximize
$$
\mathbf{Q}_{q,m}\left(N=n\right)=u_m(n)q^{n-m}(q)_m^2,
$$
i.e., as the saddle-point of $q^{m-n}(q)_m^{-2}$.  Following, for example, \cite{B20,F,P,R,V} we set $q=q(n)=e^{-\frac{1}{B\sqrt{n}}}$, which importantly is independent of $m$.  We use this $q$ throughout when referring to unimodal sequences, and we always write $m$ in terms of $r \in \mathbb{R}$ as $m=B\sqrt{n}\log(2B\sqrt{n})+Br\sqrt{n}$, and we switch freely between these notations.  Here and throughout, $r$ is assumed to be in $\frac{1}{B\sqrt{n}}(\mathbb{Z}-\log(2B\sqrt{n}) )$, so that $m \in \mathbb{Z}$. When referring to strongly unimodal sequences, we also always use\footnote{It is a bit surprising that this is also Fristedt's choice of $q$ for partitions.} $q=q(n)=e^{-\frac{1}{A\sqrt{n}}}$ here for ease of notation, and we take $r \in \frac{1}{A\sqrt{n}}(\mathbb{Z} -\log(2A\sqrt{n}))$ depending on $m$ now as $m+1=A\sqrt{n}\log(2A\sqrt{n})+Ar\sqrt{n}$. For any $[r_1,r_2] \subset \mathbb{R}$, we can break up $\mathbf{P}_{n}$ into the ranges
$$
\mathbf{P}_{n}= \left(\sum_{r < r_1} + \sum_{r \in [r_1,r_2]} + \sum_{r > r_2} \right)   \mathbf{P}_{n} \left({\rm PK}=m\right) \cdot \mathbf{P}_{n,m},
$$
and we can bound the tail ranges for any measurable set $S$ as
\begin{align}\label{E:lowertail}
\sum_{r < r_1} \mathbf{P}_{n}\left({\rm PK}=m \right) \cdot \mathbf{P}_{n,m}(S) &\leq \sum_{r < r_1} \mathbf{P}_{n}\left({\rm PK}=m\right)= \mathbf{P}_{n}\left(\frac{{\rm PK}-B\sqrt{n}\log(2B\sqrt{n})}{B\sqrt{n}}<r_1 \right),\\
\label{E:uppertail}
\sum_{r > r_2}  \mathbf{P}_{n}\left({\rm PK}=m \right) \cdot \mathbf{P}_{n,m}(S) &\leq \mathbf{P}_{n}\left(\frac{{\rm PK}-B\sqrt{n}\log(2B\sqrt{n})}{B\sqrt{n}}>r_2 \right).
\end{align}
For sequences $a_n\leq b_n$ of positive integers, we define
$$
\bm{X}_{[a_n,b_n]}:=\left(X_k^{[j]}\right)_{k \in [a_n,b_n], j \in \{L,R\}}.
$$To prove our main theorems, we show that \eqref{E:lowertail} and \eqref{E:uppertail} tend to 0 as $r_1 \to -\infty$ and $r_2 \to \infty$, respectively, and that $d_{\textrm{TV}}(\mathbf{P}_{n,m}(\bm{X}_{[a_n,b_n]}^{-1}),\mathbf{Q}_{q,m}(\bm{X}_{[a_n,b_n]}^{-1}))\to 0$ uniformly for $r \in [r_1,r_2]$ given a simple condition on $a_n$ and $b_n$.  To aid the proof of the latter, we have the following analogue of Lemma 4.2 in \cite{F}; it is proved in exactly the same way.

\begin{proposition}\label{P:B_{n,m}prop}
Let $a_n\leq b_n$ be sequences of integers and suppose that there exist $B_{n,m} \subset \mathbb{R}^{2(b_n-a_n+1)}$ such that, uniformly for $r$ in any $[r_1,r_2]$,
\begin{enumerate}[leftmargin=*,label=\textnormal{(\arabic*)}]
\item $\mathbf{Q}_{n,m}\left(\bm{X}_{[a_n,b_n]} \in B_{n,m}\right) \to 1$,
\item $\frac{\mathbf{Q}_{q,m}\left(  N=n |\bm{X}_{[a_n,b_n]}={\bm x} \right)}{\mathbf{Q}_{q,m}(N=n)} \to 1$ uniformly in ${\bm x} \in B_{n,m}$.
\end{enumerate}
Then $d_{\textnormal{TV}}(\mathbf{P}_{n,m}(\bm{X}_{[a_n,b_n]}^{-1}),\mathbf{Q}_{q,m}(\bm{X}_{[a_n,b_n]}^{-1}))\to 0$ uniformly for $r$ in $[r_1,r_2]$. For strongly unimodal sequences, the same holds if we replace $\mathbf{P}_{n,m} \mapsto \mathbf{P}^*_{n,m}$ and $\mathbf{Q}_{q,m} \mapsto \mathbf{Q}_{q,m}^*$.
\end{proposition}

So long as the sequences $a_n$ and $b_n$ fulfill a simple condition,  the sets $B_{n,m}$ above exist, as we show in Proposition \ref{P:W_ncondition}.  But first we find the asymptotic behavior of the denominator in condition (2) of Proposition \ref{P:B_{n,m}prop}. 

\begin{proposition}\label{P:denomlocalasymp}
Let $q=e^{-\frac{1}{B\sqrt{n}}}$.  Uniformly for $r$ in any $[r_1,r_2]$, we have 
\begin{align}\label{E:denominator}
\mathbf{Q}_{q,m}(N=n) &\sim \frac{1}{2 \cdot 3^{\frac 14} n^{\frac 34}},\\
\label{E:PKlocalasymp}
\mathbf{P}_{n}({\rm PK}=m) &\sim \frac{1}{B\sqrt{n}}e^{-r-e^{-r}}.
\end{align}

For strongly unimodal sequences, we have the corresponding
\begin{align}\label{E:*denominator}
\mathbf{Q}_{q,m}^*(N=n) &\sim \frac{1}{2 \cdot 6^{\frac 14} n^{\frac 34}},\\
\label{E:*PKlocalasymp}
\mathbf{P}_{n}^*({\rm PK}=m+1) &\sim \frac{1}{A\sqrt{n}}e^{-r-e^{-r}}.
\end{align}
\end{proposition}

Before proving Proposition \ref{P:denomlocalasymp}, we show that the tail ranges \eqref{E:lowertail} and \eqref{E:uppertail} tend to 0 as a corollary and we prove Theorem \ref{T:PKdist}.
\begin{corollary}\label{C:tails}
For any $[r_1,r_2] \subset \mathbb{R}$, we have
$$
\mathbf{P}_{n}\left(r_1\leq \frac{{\rm PK}-B\sqrt{n}\log(2B\sqrt{n})}{B\sqrt{n}} \leq r_2 \right) \sim e^{-e^{-r_2}}-e^{-e^{-r_1}},
$$
and thus, because $\mathbf{P}_{n}$ is a probability measure, \eqref{E:lowertail} and \eqref{E:uppertail} tend to 0 as $r_1 \to -\infty$ and $r_2 \to \infty$, respectively.  Furthermore, Theorem \ref{T:PKdist} holds.

For strongly unimodal sequences, the above holds with $\mathbf{P}_{n} \mapsto \mathbf{P}_{n}^*$ and $B \mapsto A$.
\end{corollary}

\begin{proof}
Using \eqref{E:PKlocalasymp}, we have,
\begin{multline*}
\mathbf{P}_{n}\left(\lambda \in \mathcal{U}(n) : r_1 \leq \frac{{\rm PK}(\lambda)-B\sqrt{n}\log(2B\sqrt{n})}{B\sqrt{n}} \leq r_2 \right)\\
\sim \frac{1}{B\sqrt{n}}\sum_{r \in [r_1,r_2] \cap \left( \frac{1}{B\sqrt{n}}\Z - \log(2B\sqrt{n})\right)}\hspace{-1cm} e^{-r-e^{-r}}
= \int_{r_1}^{r_2} e^{-r-e^{-r}}dr + O\left(n^{-\frac 12}\right)
=e^{-e^{-r_2}}-e^{-e^{-r_1}} +O\left(n^{-\frac 12}\right),
\end{multline*}
where the first equality comes from recognizing a Riemann sum.  Thus, the first part of Theorem \ref{T:PKdist} holds.  The second part is a consequence of the following similar calculation. Let
$$
\widetilde{{\rm PK}}:=\frac{{\rm PK}-B\sqrt{n}\log(2B\sqrt{n})}{B\sqrt{n}},
$$
so ${\rm PK}=m$ if and only if $\widetilde{{\rm PK}}=r$, and the second part of Theorem \ref{T:PKdist} is equivalent to $\mathbf{E}_n(\widetilde{{\rm PK}}) \sim \gamma.$  For this we write

\begin{align*}
\mathbf{E}_n\left(\widetilde{{\rm PK}}\right) &= \left(\sum_{r < r_1} + \sum_{ r \in [r_1,r_2]} + \sum_{r > r_2}\right) r\mathbf{P}_{n}\left(\widetilde{{\rm PK}}=r\right).
\end{align*}
where above and below the sums on $r$ run over $(\frac{1}{B\sqrt{n}}\mathbb{Z}-\log(2B\sqrt{n}) )$.  Assuming \eqref{E:PKlocalasymp}, the middle sum is
\begin{equation}\label{E:expmiddle}
\sum_{r \in [r_1,r_2]} r\mathbf{P}_{n}\left(\widetilde{{\rm PK}}=r\right)=\frac{1}{B\sqrt{n}}\sum_{r \in [r_1,r_2]} re^{-r-e^{-r}}= \int_{r_1}^{r_2} re^{-r-e^{-r}}dr + O\left(n^{-\frac 12}\right).
\end{equation}
For the upper tail, we write $m_2$ as the integer corresponding to $r_2$, and use the elementary identity $$\sum_{m > m_2} m \mathbf{P}_n({\rm PK}=m)= \sum_{m > m_2} \mathbf{P}_n({\rm PK} \geq m)+ m_2\mathbf{P}_n({\rm PK} > m_2),$$ so that
\begin{align}
\sum_{r > r_2} r\mathbf{P}_{n}\left(\widetilde{{\rm PK}}=r\right) &=\sum_{m > m_2} \frac{m-B\sqrt{n}\log\left(2B\sqrt{n}\right)}{B\sqrt{n}} \mathbf{P}_{n}\left({\rm PK}=m\right) \nonumber \\
&=\frac{1}{B\sqrt{n}} \left(\sum_{m > m_2} \mathbf{P}_n({\rm PK} \geq m)+ m_2\mathbf{P}_n({\rm PK} > m_2)\right)+\log\left(2B\sqrt{n}\right)\mathbf{P}_n({\rm PK} > m_2) \nonumber \\
&=\frac{1}{B\sqrt{n}} \sum_{m > m_2} \mathbf{P}_n\left({\rm PK} \geq m\right)+ r_2\mathbf{P}_n({\rm PK} > m_2) \nonumber \\
&=\frac{1}{B\sqrt{n}} \sum_{r > r_2} \mathbf{P}_n\left(\widetilde{{\rm PK}} \geq r\right)+ r_2\mathbf{P}_n\left(\widetilde{{\rm PK}} > r_2\right) \nonumber \\
&=\int_{r_2}^{\infty} \left(1-e^{-e^{-r}}\right)dr+r_2\left(1-e^{-e^{-r_2}}\right)+O\left(n^{-\frac{1}{2}}\right), \label{E:expuppertail}
\end{align}
where in the last step we use the first part of Theorem \ref{T:PKdist}.

In a similar way, we estimate the lower tail as
\begin{equation}\label{E:explowertail}
\sum_{r < r_1} r\mathbf{P}_{n}\left(\widetilde{{\rm PK}}=r\right)= \int_{-\infty}^{r_1}e^{-e^{-r}}dr+r_1e^{-e^{-r_1}}+O\left(n^{-\frac{1}{2}}\right).
\end{equation}
 Taking $r_1 \to -\infty$ and $r_2 \to \infty$ in \eqref{E:expmiddle}, \eqref{E:expuppertail} and \eqref{E:explowertail} gives
$$
\mathbf{E}_n\left(\widetilde{{\rm PK}}\right)\sim\int_{-\infty}^{\infty} re^{-r-e^{-r}}dr = -\int_0^{\infty} \log(t)e^{-t}dt=-\Gamma'(1)=\gamma,
$$
which proves the second part of Theorem \ref{T:PKdist}.

Repeating this entire proof with $B \mapsto A$ proves the claims for strongly unimodal sequences.
\end{proof}

We are now ready to prove Proposition \ref{P:denomlocalasymp}.

\begin{proof}[Proof of Proposition \ref{P:denomlocalasymp}] 
First note that
\begin{equation}\label{E:QprobN=n}
\mathbf{Q}_{q,m}(N=n)= (q)_m^2 q^{n-m} u_m(n).
\end{equation}
To estimate the right-hand side, we could take the probabilistic approach given in the proofs of Proposition 4.5 in \cite{F} and Proposition 3 in \cite{R}, and the reader is invited to examine these for a more heuristic approach.  Instead, we use Proposition \ref{P:saddlepointmethod}, a version of the classical saddle-point method, after representing $u_m(n)$ as a Cauchy contour integral, but this approach is essentially equivalent to Fristedt's probabilistic approach.  We omit some details when they are very similar to calculations carried out in \cite{B20} and \cite{R}.
First, write
$$
u_m(n)=\frac{1}{2\pi i} \int_{\mathcal{C}} \frac{\zeta^{m-n-1}}{(\zeta)_m^2}d\zeta,
$$
where $\mathcal{C}$ is a circle centered at 0 with radius less than 1 orientated counterclockwise.  Substituting $\zeta=e^{-\frac{1}{B\sqrt{n}}+2 \pi i \theta}$, we have
\begin{equation*}\label{E:Cauchyint}
u_m(n)=\int_{-\frac{1}{2}}^{\frac{1}{2}} \frac{e^{-\frac{m-n}{B\sqrt{n}}+ 2 \pi i(m-n)  \theta}}{\left(e^{-\frac{1}{B\sqrt{n}}+2 \pi i \theta}\right)_m^2}d \theta = \int_{-\frac{1}{2}}^{\frac{1}{2}} \exp\left(f(2 \pi i \theta)\right) d \theta,
\end{equation*}
where
\begin{equation*}\label{E:fdef}
f(z):= \frac{n-m}{B\sqrt{n}}+ (m-n) z -2 \sum_{k=1}^m \Log\left(1-e^{-\frac{k}{B\sqrt{n}}+ k z  } \right).
\end{equation*}
For ease of notation, we omit the dependence of $f$ on $n$ throughout the proof.  To prove that condition (i) in Proposition \ref{P:saddlepointmethod} holds, we need to find the asymptotic behavior of $f(0), f'(0)$, and $f''(0)$.  These are easily computed:
\begin{align*}
f(0) &= \frac{n-m}{B\sqrt{n}}-2\sum_{k=1}^m \log\left(1-e^{-\frac{k}{B\sqrt{n}}}\right), \qquad f'(0) = m-n+2\sum_{k=1}^m \frac{k}{e^{\frac{k}{B\sqrt{n}}}-1},\\
f''(0) &= 2\sum_{k=1}^m \frac{k^2 e^{-\frac{k}{B\sqrt{n}}}}{\left(1-e^{-\frac{k}{B\sqrt{n}}}\right)^2}.
\end{align*}

Next, we use Lemma \ref{L:Romik} with $x\mapsto e^{-\frac{k}{B\sqrt{n}}}$, $s \mapsto 2\pi k\th$ to compute
\begin{align*}
&\left|f(2\pi i \theta) - f(0)-  f'(0) 2 \pi i \theta -  f''(0) \frac{(2 \pi i \theta)^2}{2}\right|\\
&\hspace{1cm}=\left|\frac{n-m}{B\sqrt{n}} + (m-n) 2\pi i\th - 2\sum_{k=1}^m \log\left(1-e^{-\frac{k}{B\sqrt{n}}+k\th}\right) - \frac{n-m}{B\sqrt{n}} - 2\sum_{k=1}^m \log\left(1-e^{-\frac{k}{B\sqrt{n}}}\right)\right.\\
&\hspace{4cm}\left. -\left((m-n) + 2\sum_{k=1}^m \frac{k}{e^{\frac{k}{B\sqrt{n}}}-1}\right)2\pi i\th - \sum_{k=1}^m \frac{k^2 e^{-\frac{k}{B\sqrt{n}}}}{\left(1-e^{\frac{k}{B\sqrt{n}}}\right)^2}(2\pi i\th)^2\right|\\
&\hspace{1cm}= 2\left|{\vphantom{\sum_{k=1}^m \left(\Log\left(\frac{1-e^{-\frac{k}{B\sqrt{n}}}}{1-e^{-\frac{k}{B\sqrt{n}}+2 \pi i k \theta}}\right)- 2 \pi i k\theta \frac{e^{-\frac{k}{B\sqrt{n}}}}{1-e^{-\frac{k}{B\sqrt{n}}}} +\frac{(2\pi i k \theta)^2}{2} \frac{e^{-\frac{k}{B\sqrt{n}}}}{\left(1-e^{-\frac{k}{B\sqrt{n}}}\right)^2}\right)}}\right.
\sum_{k=1}^m
\left({\vphantom{\Log\left(\frac{1-e^{-\frac{k}{B\sqrt{n}}}}{1-e^{-\frac{k}{B\sqrt{n}}+2 \pi i k \theta}}\right)- 2 \pi i k\theta \frac{e^{-\frac{k}{B\sqrt{n}}}}{1-e^{-\frac{k}{B\sqrt{n}}}} +\frac{(2\pi i k \theta)^2}{2} \frac{e^{-\frac{k}{B\sqrt{n}}}}{\left(1-e^{-\frac{k}{B\sqrt{n}}}\right)^2}}}\right.
\Log\left(\frac{1-e^{-\frac{k}{B\sqrt{n}}}}{1-e^{-\frac{k}{B\sqrt{n}}+2 \pi i k \theta}}\right)-  \frac{e^{-\frac{k}{B\sqrt{n}}}}{1-e^{-\frac{k}{B\sqrt{n}}}}2 \pi i k\theta + \frac{e^{-\frac{k}{B\sqrt{n}}}}{\left(1-e^{-\frac{k}{B\sqrt{n}}}\right)^2}
\left.{\vphantom{\Log\left(\frac{1-e^{-\frac{k}{B\sqrt{n}}}}{1-e^{-\frac{k}{B\sqrt{n}}+2 \pi i k \theta}}\right)- 2 \pi i k\theta \frac{e^{-\frac{k}{B\sqrt{n}}}}{1-e^{-\frac{k}{B\sqrt{n}}}} +\frac{(2\pi i k \theta)^2}{2} \frac{e^{-\frac{k}{B\sqrt{n}}}}{\left(1-e^{-\frac{k}{B\sqrt{n}}}\right)^2}}}\frac{(2\pi i k \theta)^2}{2}\right)
\left.{\vphantom{\sum_{k=1}^m 
\left(\Log\left(\frac{1-e^{-\frac{k}{B\sqrt{n}}}}{1-e^{-\frac{k}{B\sqrt{n}}+2 \pi i k \theta}}\right)
- 2 \pi i k\theta \frac{e^{-\frac{k}{B\sqrt{n}}}}{1-e^{-\frac{k}{B\sqrt{n}}}}
+\frac{(2\pi i k \theta)^2}{2} 
\frac{e^{-\frac{k}{B\sqrt{n}}}}{\Bigg(1-e^{-\frac{k}{B\sqrt{n}}}\Bigg)^2}
\frac{(2\pi i k \theta)^2}{2}\right)}}\right|\\
&\hspace{1cm}\leq C \sum_{k=1}^m \frac{k^3e^{-\frac{k}{B\sqrt{n}}}}{\left(1-e^{-\frac{k}{B\sqrt{n}}}\right)^3} |\theta|^3,
\end{align*}
for some constant $C$.  Recognizing Riemann sums for a convergent integral, the sum is bounded by
\begin{align}
 B^4n^2 \sum_{k \geq 1}\frac{\left(\frac{k}{B\sqrt{n}}\right)^3e^{-\frac{k}{B\sqrt{n}}}}{\left(1-e^{-\frac{k}{B\sqrt{n}}}\right)^3} \frac{1}{B\sqrt{n}} &=    B^4n^2 \int_{0}^{\infty} \frac{u^3e^{-u}}{(1-e^{-u})^3}du \left(1+O\left(\frac{1}{\sqrt{n}}\right)\right) =O\left(n^2\right), \label{E:Taylorerror}
\end{align}
where the constant is independent of $n$.  Now we work out the asymptotic behavior of $f(0)$, $f'(0)$, and $f''(0)$.  The details here are very similar to the proofs of Propositions 1--3 in \cite{R}.
For $f(0)$, we take $g(u):= -\log(1-e^{-\frac{u}{B\sqrt{n}}})$ in \eqref{E:eulermac1}.  Estimating as in the proof of Proposition 1 in \cite{R} gives       
\begin{align}
f(0) = 2\pi\sqrt{\frac{n}{3}}- \log(n) -r-e^{-r}-\log\left(\frac{12}{\pi} \right) +o(1).  \label{E:f(0)asymp}
\end{align}

The asymptotic behavior of the derivatives is similar to the proof of Proposition 2 in \cite{R}.  For $f'(0)$, we take $g(u):=\frac{u}{e^{\frac{u}{B\sqrt{n}}}-1}$ with \eqref{E:eulermac2}, noting that 
$$
g'(u) = \frac{e^{\frac{u}{B\sqrt{n}}}-\frac{u}{B\sqrt{n}}e^{\frac{u}{B\sqrt{n}}}-1}{\left(e^{\frac{u}{B\sqrt{n}}}-1\right)^2}.
$$  A short calculation gives 
\begin{align}
f'(0) =O\left(\sqrt{n} \log (n) \right). \label{E:f'(0)asymp}
\end{align}

Finally, for $f''(0)$ take $g(u)= \frac{u^2e^{-\frac{u}{B\sqrt{n}}}}{(1-e^{-\frac{u}{B\sqrt{n}}})^2}$ in \eqref{E:eulermac2}.  Using that 
$$
g'(u)=-\frac{ue^{-\frac{u}{B\sqrt{n}}}\left(\left(\frac{u}{B\sqrt{n}}+2\right) e^{-\frac{u}{B\sqrt{n}}}+\frac{u}{B\sqrt{n}}-2\right)}{\left(1-e^{-\frac{u}{B\sqrt{n}}}\right)^3},
$$
a short calculation gives 
\begin{align}
f''(0)= \frac{2\sqrt{3}}{\pi}n^{\frac 32} + O\left(n\log(n)^2\right). \label{E:f''(0)asymp}
\end{align}
Thus, assumption (i) in Proposition \ref{P:saddlepointmethod} holds for any fixed $\varepsilon >0$ and $|\theta| \leq \varepsilon n^{-\frac{1}{2}}$.

It remains to check that assumption (ii) in Proposition \ref{P:saddlepointmethod} holds.  First, we write
\begin{align*}
 \text{Re} f(2\pi i \theta)-f(0) 
&=  -2\sum_{k=1}^m \left( \text{Re} \left(\Log\left(1-e^{-\frac{k}{B\sqrt{n}}+2 \pi i \theta}\right) \right) 
   - \Log\left(1-e^{-\frac{k}{B\sqrt{n}}}\right)\right)   \\
&= 2\sum_{k=1}^m \sum_{\ell \geq 1}  \frac{e^{-\ell \frac{k}{B\sqrt{n}}}}{\ell} \left( \cos(2 \pi k \ell \theta)-1 \right) 
\leq 2\sum_{k=1}^m  e^{-\frac{k}{B\sqrt{n}}} \left( \cos(2 \pi k \theta)-1 \right).
\end{align*}

We analyze the sum in the exponential as in \cite{R}, p. 13.  Because it is a Riemann sum, we have, for $n$ sufficiently large and $\varepsilon n^{-\frac 12}\leq |\theta| \leq \frac{1}{2}$,
\begin{align*}
\sum_{k=1}^m  &e^{- \frac{k}{B\sqrt{n}}} \left( \cos(2 \pi k \theta)-1 \right)\\
&=-B\sqrt{n}\sum_{k=1}^m  e^{- \frac{k}{B\sqrt{n}}} \left( 1-\cos(2 \pi k \theta) \right)  \frac{1}{B\sqrt{n}}
< -\frac{B\sqrt{n}}{2}\int_0^{\frac{m}{B\sqrt{n}}} e^{-u}\left(1-\cos\left(2 \pi B\sqrt{n}\theta u\right)\right)du\\
&< -\frac{B\sqrt{n}}{2}\inf_{\theta \geq \varepsilon n^{-\frac{1}{2}}}\int_0^{T} e^{-u}\left(1-\cos\left(2 \pi B\sqrt{n}\theta u\right)\right)du
= -\frac{B\sqrt{n}}{2}\inf_{s \geq \varepsilon}\int_0^{T} e^{-u}\left(1-\cos\left(2 \pi Bs u\right)\right)du,
\end{align*}
for any $T >0$, since $\frac{m}{B\sqrt{n}} \to \infty$.  The infimum above is positive, since the function
\[
s \mapsto \int_0^{T} e^{-u}\left(1-\cos\left(2 \pi Bs u\right)\right)du
\]
is continuous and nonzero on $[\varepsilon, \infty)$ and tends to $1-e^{-T}>0$ as $s \to \infty$ by the Riemann--Lebesgue Lemma, so assumption (ii) in Proposition \ref{P:saddlepointmethod} holds, and we conclude that
\begin{equation}\label{E:u(n,m)asymp}
u_m(n)\sim \frac{e^{f(0)}}{\sqrt{2\pi f''(0)}}.
\end{equation}
Recalling \eqref{E:QprobN=n} and using $e^{f(0)}=(q)_m^{-2} q^{m-n}$, the cancellation yields
$$
\mathbf{Q}_{q,m}(N=n)\sim \frac{1}{\sqrt{2\pi f''(0)}}.
$$
Plugging in \eqref{E:f''(0)asymp} proves \eqref{E:denominator}.

Note that the asymptotic behavior of $f(0)$ is only required to prove \eqref{E:PKlocalasymp}.  By definition,
$$
\mathbf{P}_{n}({\rm PK}=m)= \frac{u_m(n)}{u(n)}.
$$
Plugging in \eqref{E:f(0)asymp} and \eqref{E:f''(0)asymp} into \eqref{E:u(n,m)asymp} and using the asymptotic due to Auluck (\cite{A}, eq. (24)) $u(n) \sim \frac{e^{2 \sqrt{\frac n3}}}{8 \cdot 3^{\frac 34} n^{\frac 54}}$ proves \eqref{E:PKlocalasymp}. 

For strongly unimodal sequences, we make a very similar argument using Proposition \ref{P:saddlepointmethod}.  A more probabilistic approach for this case is similar to the proof of Proposition 3 in \cite{B20}.
We write
$$
u^*_m(n)= \int_{-\frac{1}{2}}^{\frac{1}{2}} \exp\left(f_*(2\pi i \theta)\right)d \theta,
$$
where this time,
$$
f_*(z):=\frac{n-m-1}{A\sqrt{n}}+(m+1-n)z + 2 \sum_{k=1}^m \log\left( 1+ e^{-\frac{k}{A\sqrt{n}} + k z }\right).
$$
We have
\begin{align*}
f_*(0) &= \frac{n-m-1}{A\sqrt{n}}+2\sum_{k=1}^m \log\left(1+e^{-\frac{k}{A\sqrt{n}}}\right), \qquad f_*'(0) = m+1-n+2\sum_{k=1}^m \frac{k}{e^{\frac{k}{A\sqrt{n}}}+1},\\
f''_*(0) &= 2\sum_{k=1}^m \frac{k^2 e^{-\frac{k}{A\sqrt{n}}}}{\left(1+e^{-\frac{k}{A\sqrt{n}}}\right)^2}.
\end{align*}

Next, we use Lemma \ref{L:Romik} with $x\mapsto e^{-\frac{k}{A\sqrt{n}}}$, $s \mapsto 2\pi k\th$ to compute
\begin{align*}
&\left|f_*(2\pi i \theta) - f_*(0)-  f_*'(0)2 \pi i \theta -  f_*''(0)\frac{(2 \pi i \theta)^2}{2} \right| \leq C \sum_{k=1}^m \frac{k^3e^{-\frac{k}{A\sqrt{n}}}}{\left(1-e^{-\frac{k}{A\sqrt{n}}}\right)^3}|\theta|^3 =O\left(\theta^3n^2\right),
\end{align*}
similarly to before.  We now find the asymptotic behavior of $f_*(0)$, $f_*'(0)$, and $f_*''(0)$ using Lemma \ref{L:eulermac}.  The reader may consult the proof of Propositions 1 and 2 in \cite{B20} for similar asymptotic calculations.  We obtain the following
\begin{align}
f_*(0)&=\pi\sqrt{\frac{2}{3}n} - r -e^{-r}-\log\left(4A\sqrt{n}\right)+o(1), \label{E:*f(0)asymp} \\
f_*'(0)& = O\left(\sqrt{n} \log (n)\right), \label{E:*f'(0)asymp} \\
f_*''(0)&= \frac{2\sqrt{6}}{\pi}n^{\frac 32} + O\left(n\log(n)^2\right). \label{E:*f''(0)asymp}
\end{align}

It remains to check condition (ii).  Here, the analysis is similar to (p. 253 in \cite{RS}).  First, we use the identity
$$
\left|\frac{1+q^ke^{2\pi i k\theta}}{1+q^k}\right|^2=\frac{1}{\left(1+q^k\right)^2}\left(1+2q^k\cos\left(2 \pi k \theta\right)+q^{2k}\right)=1-\frac{2q^k\left(1-\cos(2\pi \theta)\right)}{\left(1+q^k\right)^2},
$$
where we recall that $q=e^{-\frac{1}{A\sqrt{n}}}$.  Thus,
\begin{align*}
\left| \exp\left(f_*(2\pi i \theta)-f_*(0)\right) \right| 
&=\prod_{k=1}^m\left|   \frac{1+q^ke^{2 \pi i k \theta}}{1+q^k}\right|^2=\prod_{k=1}^m\left(1-\frac{2q^k\left(1-\cos(2\pi \theta)\right)}{\left(1+q^k\right)^2}\right),
\end{align*}
so
\begin{align*}
\mathrm{Re}(f_*(2\pi i \theta)) - f_*(0)
&= \sum_{k=1}^m\log\left(1-\frac{2q^k(1-\cos(2 \pi k \theta))}{\left(1+q^k\right)^2} \right) \leq -\sum_{k=1}^m  \frac{2q^k(1-\cos(2 \pi k \theta))}{(1+q^k)^2}  \\  
&\leq -2\sum_{k=1}^m  q^k(1-\cos(2 \pi k \theta))  ,
\end{align*}
which gives the same Riemann sum as in the estimate of $\mathrm{Re}(f(2\pi i \theta)) - f(0)$ with the constants $A$ and $B$ swapped.  Thus, the above is $-\delta_{\varepsilon}\sqrt{n}$ for some $\delta_{\varepsilon}>0$ and $\varepsilon n^{-\frac{1}{2}} \leq |\theta| \leq \frac{1}{2}$ exactly as before, so assumption (ii) holds, and Proposition \ref{P:saddlepointmethod} implies
$$
u_m^*(n) \sim \frac{e^{f_*(0)}}{\sqrt{2 \pi f_*''(0)}}.
$$
We then get \eqref{E:*denominator} exactly as before, by plugging in \eqref{E:*f''(0)asymp}.  For \eqref{E:*PKlocalasymp}, we plug in \eqref{E:*f(0)asymp} and \eqref{E:*f''(0)asymp} and use the asymptotic $u^*(n) \sim \frac{e^{\pi \sqrt{\frac{2}{3}n}}}{8 \cdot 6^{\frac 14} n^{\frac 34}}$ due to Rhoades (\cite{Rh}, Theorem 1.1).
\end{proof}

The next proposition handles the numerator in Proposition \ref{P:B_{n,m}prop} by providing an explicit set $B_{n,m}$, given a simple condition.

\begin{proposition}\label{P:W_ncondition}
Suppose that $a_n \leq b_n$ are sequences of integers such that
\begin{equation}\label{E:b_ncondition}
\sum_{a_n \leq k \leq b_n} \frac{k^2q^k}{\left(1-q^k\right)^2}=o\left(c_n^2\right)
\end{equation}
holds for a sequence $c_n=o(n^{\frac{3}{4}})$.  Then for $b_{n,m}:=\min\{b_n,m\}$,
$$
B_{n,m}:= \left\{ \left(x^{[j]}_k\right)_{\substack{j \in \{L,R\} \\ a_n \leq k \leq b_{n,m}}} : \left|\sum_{a_n \leq k \leq b_{n,m}} \frac{2kq^k}{1-q^k} - \sum_{a_n \leq k \leq b_{n,m}} k\left(x^{[L]}_k+x^{[R]}_k\right) \right| \leq c_n \right\} \times \{0\}^{2(b_n-b_{n,m})}
$$
satisfies the hypotheses of Proposition \ref{P:B_{n,m}prop}, so $d_{\textnormal{TV}}(\mathbf{P}_{n,m}(\bm{X}_{[a_n,b_n]}^{-1}),\mathbf{Q}_{q,m}(\bm{X}_{[a_n,b_n]}^{-1}))\to 0$ uniformly for $r$ in any $[r_1,r_2]$.

For strongly unimodal sequences, we have $d_{\textnormal{TV}}(\mathbf{P}^*_{n,m}(\bm{X}_{[a_n,b_n]}^{-1}),\mathbf{Q}^*_{q,m}(\bm{X}_{[a_n,b_n]}^{-1}))\to 0$ if
\begin{equation}\label{E:b_ncondition}
\sum_{a_n \leq k \leq b_n} \frac{k^2q^k}{\left(1+q^k\right)^2}=o\left(c_n^2\right)
\end{equation}
for $c_n=o(n^{\frac{3}{4}})$ with
$$
B^*_{n,m}:= \left\{ \left(x^{[j]}_k\right)_{\substack{j \in \{L,R\} \\ a_n \leq k \leq b_{n,m}}} : \left|\sum_{a_n \leq k \leq b_{n,m}} \frac{2kq^k}{1+q^k} - \sum_{a_n \leq k \leq b_{n,m}} k\left(x^{[L]}_k+x^{[R]}_k\right) \right| \leq c_n \right\} \times \{0\}^{2(b_n-b_{n,m})}.
$$
\end{proposition}

\begin{proof}
To show $\mathbf{Q}_{n,m}\left(\bm{X}_{[a_n,b_n]} \in B_{n,m}\right) \to 1$, we use Lemma \ref{L:Chebyshev}.  Note that we have
$$
\mathbf{E}_{q,m}\left(X_{k}^{[L]}\right) = \left(1-q^k\right) \sum_{\ell \geq 1} \ell q^{k\ell}=\frac{q^k}{1-q^k},
$$
$$
\mathbf{Var}_{q,m}\left(X_{k}^{[L]}\right) = \left(1-q^k\right) \sum_{\ell \geq 1} \ell^2q^{k\ell} - \frac{q^{2k}}{\left(1-q^k\right)^2}=\frac{q^k+q^{2k}}{\left(1-q^k\right)^2}- \frac{q^{2k}}{\left(1-q^k\right)^2}=\frac{q^k}{\left(1-q^k\right)^2}.
$$
These are also exactly the same for $X_{k}^{[R]}$.  Thus, we have
$$
\mathbf{E}_{q,m}\left(\sum_{a_n \leq k \leq b_{n,m}} k\left(X_{k}^{[L]}+X_{k}^{[R]}\right)\right)=\sum_{a_n \leq k \leq b_{n,m}} \frac{2kq^k}{1-q^k},
$$
and, using independence,
$$
\mathbf{Var}_{q,m}\left(\sum_{a_n \leq k \leq b_{n,m}} k\left(X_{k}^{[L]}+X_{k}^{[R]}\right)\right)=\sum_{a_n \leq k \leq b_{n,m}} \frac{2k^2q^k}{(1-q^k)^2}.
$$
Thus, by Chebyshev's Inequality (Lemma \ref{L:Chebyshev}) and the definition of $B_{n,m}$
$$
\mathbf{Q}_{n,m}\left(\bm{X}_{[a_n,b_n]} \in \mathbb{R}^{b_n-a_n+1} \setminus B_{n,m}\right) \leq c_n^{-2}\sum_{a_n \leq k \leq b_{n,m}} \frac{2k^2q^k}{\left(1-q^k\right)^2} = o(1),
$$
as required. Now let ${\bm x}=(x^{[j]}_k)_{j \in \{L,R\}, a_n \leq k \leq b_{n,m}} \times \{0\}^{b_n-b_{n,m}} \in B_{n,m}$, and write 
$$
\sum {\bm x}:=\sum_{a_n \leq k \leq b_{n,m}} k\left(x^{[L]}_k+x^{[R]}_k\right).
$$
We want to show that
$$
\mathbf{Q}_{q,m}(N=n | \bm{X}_{[a_n,b_n]}={\bm x}) \sim \frac{1}{2 \cdot 3^{\frac 14}n^{\frac 34}}.
$$
The proof of this is simply a slight adjustment to the proof of the first part of Proposition \ref{P:denomlocalasymp}.  Again, we wish to apply Proposition \ref{P:saddlepointmethod} to a Cauchy integral.  We write, using independence,
\begin{align*}
\mathbf{Q}_{q,m}(N=n | \bm{X}_{[a_n,b_n]}={\bm x})&= \frac{\mathbf{Q}_{q,m}\left(N=n \ \text{and} \ \bm{X}_{[a_n,b_n]}={\bm x} \right)}{\mathbf{Q}_{q,m}(\bm{X}_{[a_n,b_n]}={\bm x})} \\ 
&= \frac{\mathbf{Q}_{q,m}\left(N=n \ \text{and} \ \bm{X}_{[a_n,b_n]}={\bm x} \right)}{\prod_{a_n \leq k \leq b_{n,m}} \mathbf{Q}_{q,m}\left(X_k^{[L]}=x^{[L]}_k\right) \mathbf{Q}_{q,m}\left(X_k^{[R]}=x^{[R]}_k\right)} \\ 
&= \frac{(q)_m^2q^{n-m}}{q^{\sum {\bm x}}\frac{(q)^2_{b_{n,m}}}{(q)^2_{a_n-1}}} \#\left\{ \lambda : N(\lambda)=n, {\rm PK}(\lambda)=m, \bm{X}_{[a_n,b_n]}(\lambda)={\bm x} \right\} \\
&= \frac{(q)_m^2(q)^2_{a_n-1}q^{n-m}}{(q)^2_{b_{n,m}}q^{\sum {\bm x}}} \int_{-\frac 12}^{\frac 12} \exp(F(2 \pi i \theta))d \theta,
\end{align*}
where now 
$$
F(z):=\frac{n-m-\sum {\bm x}}{B\sqrt{n}}+ \left(\sum {\bm x} + m-n\right) z -2 \sum_{k \in [1,m] \setminus [a_n,b_{n,m}]} \Log\left(1-e^{-\frac{k}{B\sqrt{n}}+ k z}\right).
$$
We have, with all sums are over $k \in [1,m] \setminus [a_n,b_{n,m}]$
\begin{align*}
F(0) &= \frac{n-m-\sum {\bm x}}{B\sqrt{n}}-2\sum_k \log\left(1-e^{-\frac{k}{B\sqrt{n}}}\right), \qquad F'(0) = \sum {\bm x} + m-n+2\sum_k \frac{k}{e^{\frac{k}{B\sqrt{n}}}-1},\\
F''(0) &= 2\sum_k \frac{k^2 e^{-\frac{k}{B\sqrt{n}}}}{\left(1-e^{-\frac{k}{B\sqrt{n}}}\right)^2},
\end{align*}
The exact asymptotic behavior of $F(0)$ is immaterial, since
$$
\frac{(q)_m^2(q)^2_{a_n-1}q^{n-m}}{(q)^2_{b_{n,m}}q^{\sum {\bm x}}}e^{F(0)}=1.
$$
Since ${\bm x} \in B_{n,m}$, we have
$$
F'(0) = m-n +2\sum_{k \leq m} \frac{k}{e^{\frac{k}{B\sqrt{n}}}-1} + o(c_n)=f'(0)+o(c_n) \sim f'(0),
$$
and by hypothesis,
$$
F''(0)=2\sum_{k\leq m} \frac{k^2 e^{-\frac{k}{B\sqrt{n}}}}{\left(1-e^{-\frac{k}{B\sqrt{n}}}\right)^2} +o\left(c_n^2\right)=f''(0)+o\left(c_n^2\right) \sim f''(0).
$$
Finally, Lemma \ref{L:Romik} is used on individual summands as in the proof of Proposition \ref{P:denomlocalasymp}, so assumption (i) in Proposition \ref{P:saddlepointmethod} holds.  

It remains to check assumption (ii).  First, we claim that there is an interval $[\alpha\sqrt{n},\beta\sqrt{n}]$ contained in $[1,m] \setminus [a_n,b_{n,m}]$.  Indeed, it is contained in $[1,m]$ for large $n$ and $r$ in any $[r_1,r_2]$, and we have
$$
\sum_{\alpha\sqrt{n} \leq k \leq \beta\sqrt{n}} \frac{k^2q^k}{\left(1-q^k\right)^2} \asymp n^{\frac 32},
$$
by estimating the Riemann sum by an integral as in the calculation for $f''(0)$.  If for all $\alpha$ and $\beta$ this intersected $[a_n,b_{n,m}]$, then this would contradict the condition \eqref{E:b_ncondition}.

Now we can estimate the contribution from the minor arc as before, using
$$
\sum_{k \in [1,m] \setminus [a_n,b_{n,m}]} e^{-\frac{k}{B\sqrt{n}}}\left(\cos\left(2 \pi k \theta \right)-1 \right)<\sum_{\alpha\sqrt{n} \leq k \leq \beta\sqrt{n}}e^{-\frac{k}{B\sqrt{n}}}\left(\cos\left(2 \pi k \theta \right)-1 \right).
$$
Similarly to before, we bound this negative sum from above by (for $\varepsilon n^{-\frac{1}{2}} \leq |\theta| \leq \frac{1}{2}$)
$$
-\frac{\sqrt{n}}{2B}\inf_{s \geq \varepsilon}\int_{\frac{\alpha}{B}}^{\frac{\beta}{B}}\left(1-\cos(2 \pi s u) \right) e^{-u} du,
$$
 and assumption (ii) holds as before.  Employing Proposition \ref{P:saddlepointmethod}, we conclude that for ${\bm x} \in B_{n,m}$,
$$
\mathbf{Q}_{q,m}(N=n | \bm{X}_{[a_n,b_n]}={\bm x} )\sim \frac{(q)_m^2(q)^2_{a_n-1}q^{n-m}}{(q)^2_{b_{n,m}}q^{\sum {\bm x}}}e^{F(0)} \frac{1}{\sqrt{2\pi F''(0)}} \sim \frac{1}{\sqrt{2\pi f''(0)}} \sim \mathbf{Q}_{q,m}(N=n).
$$

The proof for strongly unimodal sequences is similar, mirroring the proof of Proposition \ref{P:denomlocalasymp} for strongly unimodal sequences.  The main difference is that, if $\mathbf{E}^*_{q,m}$ and $\mathbf{Var}^*_{q,m}$ denote expectation and variance under $\mathbf{Q}_{q,m}^*$, then as $X_{k}^{[j]} \in \{0,1\}$, we have
$$
\mathbf{E}^*_{q,m}\left(X_{k}^{[j]}\right)=\frac{q^k}{1+q^k} \qquad \text{and} \qquad \mathbf{Var}^*_{q,m}\left(X_{k}^{[j]}\right)=\frac{q^k}{(1+q^k)^2},
$$
for $j \in \{L,R\}$.
\end{proof}

\section{Proofs of the main results}\label{S:details}

We now apply Proposition \ref{P:W_ncondition} to the sequences $a_n$ and $b_n$ for which $\bm{X}_{[a_n,b_n]}$ determines the random variables in Theorems \ref{T:largeparts}, \ref{T:smallparts}, \ref{T:totalsmallparts}, \ref{T:*largesmallparts}, and \ref{T:*smallparts}.  We then use the conditioned Boltzmann model $\mathbf{Q}_{q,m}$ to compute probability densities, identify Riemann sums and conclude the limiting distributions.

\subsection{Small parts: the proofs of Theorem \ref{T:smallparts}, Corollary \ref{C:smallpartsskewness}, and Theorem \ref{T:*smallparts}}\label{S:smallparts}

\begin{proof}[Proof of Theorem \ref{T:smallparts}]  Take $a_n=1$ and $b_n=k_n=o(n^{\frac{1}{4}})$, so that 
$$
\sum_{k \leq k_n} \frac{k^2q^k}{\left(1-q^k\right)^2} \sim \frac{n^{\frac 32}}{B} \int_0^{\frac{k_n}{B\sqrt{n}}} \frac{u^2e^{-u}}{(1-e^{-u})^2}du=o\left(n^{\frac 32}\right),
$$
and Proposition \ref{P:W_ncondition} applies to $\bm{X}_{[1,k_n]}$. Let ${\bm x}=(x^{[j]}_k)_{j \in \{L,R\}, k \leq k_n} \in \mathbb{R}^{2k_n}_{\geq 0}$ be such that $x^{[j]}_k\in \frac{k}{B\sqrt{n}}\mathbb{N}.$  Then, using independence and the exact same analysis as in \cite{F} \S 5,
\begin{align}
\mathbf{Q}_{q,m}\left(\bm{X}_{[1,k_n]} = {\bm x}\right)&=\prod_{\substack{j \in \{L,R\} \\ k \leq k_n}} \mathbf{Q}_{q,m}\left(X_{k}^{[j]} = \frac{B\sqrt{n}}{k}x^{[j]}_k\right)
= \prod_{\substack{j \in \{L,R\} \\ 1 \leq k \leq k_n}} \left(1-q^k\right)q^{\frac{B\sqrt{n}}{k}x^{[j]}_kk} \nonumber \\
&\sim  k_n!^2 \left(\frac{1}{B\sqrt{n}}\right)^{2k_n}\prod_{\substack{j \in \{L,R\} \\ k \leq k_n}} e^{-x^{[j]}_k}
=\prod_{\substack{j \in \{L,R\} \\ k \leq k_n}} e^{-x^{[j]}_k}  \frac{k}{B\sqrt{n}}. \label{E:smpartsQprobs}
\end{align}
Note that the second to last step, which comes from the analysis in \cite{F}, is the only place $k_n=o(n^{\frac 14})$ is used, rather than just $o(n^{\frac 12})$.  This is uniform in ${\bm x}$ and independent of $m$.  Hence, for any
$$
B=\prod_{\substack{k \leq k_n \\ j \in \{L,R\}}} \left.\left(-\infty,v_{k}^{[j]}\right.\right] \subset \mathbb{R}^{2k_n},
$$we have the following, uniformly for $r$ in any $[r_1,r_2]$, by recalling that $x^{[j]}_k \in \frac{k}{B\sqrt{n}}\mathbb{N}$:
\begin{align*}	
\mathbf{P}_{n,m}\left(\bm{X}_{[1,k_n]} \in B\right) &\sim \mathbf{Q}_{q,m}\left(\bm{X}_{[1,k_n]} \in B\right)\sim \sum_{{\bm w} \in B}\prod_{\substack{j \in \{L,R\} \\ k \leq k_n}} e^{-w^{[j]}_k} \frac{k}{B\sqrt{n}}\\ &\sim \prod_{\substack{k \leq k_n \\ j \in \{L,R\}}} \int_{-\infty}^{v_{k}^{[j]}} e^{- u_{k}^{[j]}} du_{k}^{[j]}=:\bm{\nu}_n(B).
\end{align*}
Thus by Corollary \ref{C:tails}, we have
\begin{align*}
&|\mathbf{P}_{n}\left(\bm{X}_{[1,k_n]} \in B\right)-\bm{\nu}_n(B)|  \\ &=\left| \sum_{r} \mathbf{P}_{n,m}\left(\bm{X}_{[1,k_n]} \in B\right)\mathbf{P}_{n}({\rm PK}=m) - \bm{\nu}_n(B)\right| \\
&=\left| \sum_{r} \left(\mathbf{P}_{n,m}\left(\bm{X}_{[1,k_n]} \in B\right)- \bm{\nu}_n(B) \right)\mathbf{P}_{n}({\rm PK}=m) \right| \\
&\leq \left(\sum_{r < r_1} + \sum_{r \in [r_1,r_2]} + \sum_{r> r_2} \right) \left|\mathbf{P}_{n,m}\left(\bm{X}_{[1,k_n]} \in B\right)-\bm{\nu}_n(B)\right|\mathbf{P}_{n}({\rm PK}=m) \\
&\leq 2e^{-e^{-r_1}} +  \sum_{r \in [r_1,r_2]}\left|\mathbf{P}_{n,m}\left(\bm{X}_{[1,k_n]} \in B\right)-\bm{\nu}_n(B)\right|\mathbf{P}_{n}({\rm PK}=m) + 2\left(1-e^{-e^{-r_2}}\right) \\
&\sim  2e^{-e^{-r_1}} +2\left(1-e^{-e^{-r_2}}\right).
\end{align*}
Taking $r_1$ and $r_2$ arbitrarily close to $-\infty$ and $\infty$, respectively, completes the proof of the first part of Theorem \ref{T:smallparts}. 
The second part is proved in the same way, merely noting that for $k=o(n^{\frac{1}{2}})$ and $w^{[L]},w^{[R]} \in \frac{k}{B\sqrt{n}}\mathbb{N}_0$, we have
\begin{equation*}
\mathbf{Q}_{q,m}\left(\frac{kX_k^{[L]}}{B\sqrt{n}}=w^{[L]}, \frac{kX_k^{[R]}}{B\sqrt{n}}=w^{[R]}\right)=\left(1-e^{-\frac{k}{B\sqrt{n}}}\right)^2e^{-w^{[L]}-w^{[R]}} \sim \left(\frac{k}{B\sqrt{n}}\right)^2e^{-w^{[L]}-w^{[R]}},
\end{equation*}
and similarly for $k=\lfloor c\sqrt{n} \rfloor $ and $u^{[L]}, u^{[R]} \in \mathbb{N}_0$,
\begin{equation*}
\mathbf{Q}_{q,m}\left(X_k^{[L]}=u^{[L]}, X_k^{[R]}=u^{[R]}\right)\sim \left(1-e^{-\frac{c}{B}}\right)^2e^{-\frac{cu^{[L]}}{B}-\frac{cu^{[R]}}{B}}. \qedhere
\end{equation*}
\end{proof}

Corollary \ref{C:smallpartsskewness} is proved similarly as above, again using \eqref{E:smpartsQprobs} to estimate the $\mathbf{Q}_{q,m}$ probabilities.

We now turn to small parts in strongly unimodal sequences.

\begin{proof}[Proof of Theorem \ref{T:*smallparts}]
To prove Theorem \ref{T:*smallparts}, we can apply Proposition \ref{P:W_ncondition} when we take $a_n=1$ and $b_n=k_n=o(n^{\frac{1}{2}})$, for
$$
\sum_{k \leq k_n} \frac{k^2q^k}{\left(1+q^k\right)^2} \sim \frac{n^{\frac 32}}{A} \int_0^{\frac{k_n}{A\sqrt{n}}} \frac{u^2e^{-u}}{(1+e^{-u})^2}du=o\left(n^{\frac 32}\right).
$$
Now, if ${\bm x} \in \{0,1\}^{2k_n}$, then we simply note that
\begin{align*}
\mathbf{Q}_{q,m}^*\left(\bm{X}_{[1,k_n]}={\bm x}\right) = \prod_{\substack{k \leq k_n \\ j \in \{L,R\}}} \mathbf{Q}_{q,m}^*\left(X_{k}^{[j]}=x^{[j]}_k\right)=\prod_{\substack{k \leq k_n \\ j \in \{L,R\}}} \frac{e^{-\frac{kx^{[j]}_k}{A\sqrt{n}}}}{1+e^{-\frac{k}{A\sqrt{n}}}} \sim \frac{1}{2^{2k_n}}. \tag* \qedhere
\end{align*}
 \end{proof}

\subsection{Large parts: the proofs of Theorems \ref{T:largeparts} and \ref{T:*largesmallparts}}
Let $t_n=o(n^{\frac{1}{4}})$ and let
$$
W_n:=\left(\frac{Y_{t}^{[j]}-B\sqrt{n}\log(2B\sqrt{n})}{B\sqrt{n}}\right)_{j \in \{L,R\}, 1 \leq t \leq t_n}.
$$

Define $\bm{\xi}_{n,m}:=\mathbf{P}_{n,m}(W_n^{-1})$, $\bm{\zeta}_{n,m}:=\mathbf{Q}_{q,m}(W_n^{-1})$ and let $\bm{\nu}_{n,r}$ be the probability measure on $(-\infty,r]^{2t_n}$ with density
$$
\begin{cases} \frac{1}{2^{2t_n}}e^{e^{-r}-\sum_{t=1}^{t_n} (u_t^{[L]}+u_t^{[R]}) - \frac{e^{-u_{t_n}^{[L]}}}{2}- \frac{e^{-u_{t_n}^{[R]}}}{2}}  & \text{if $u_1^{[j]} \geq \dots \geq u_{t_n}^{[j]}$ for $j \in \{L,R\}$,} \\ 0 & \text{otherwise.} \end{cases}
$$

Here it is easier to first that $\bm{\zeta}_{n,m}(U) \sim \bm{\nu}_{n,r}(U)$ uniformly for $r$ in any $[r_1,r_2]$ for
\begin{equation}\label{E:openUdef}
U=\prod_{\substack{t \leq t_n \\ j \in \{L,R\}}} \left.\left(-\infty,v_{t}^{[j]}\right.\right] .
\end{equation}
We then use Proposition \ref{P:W_ncondition} to complete the proof. Let ${\bm w}=(w^{[j]}_t)_{j \in \{L,R\}, t \leq t_n}$ be such that
$$
r\geq w^{[L]}_1 \geq \dots \geq w^{[L]}_{t_n} \qquad r\geq w^{[R]}_1 \geq \dots \geq w^{[R]}_{t_n}
$$
and 
$$
y_{t}^{[j]}:=B\sqrt{n}\left(w^{[j]}_t+\log\left(2B\sqrt{n}\right)\right) \in \mathbb{Z}.
$$
Directly we find
\begin{align*}
	\bm{\zeta}_{n,m}({\bm w}) &= \mathbf{Q}_{q,m} \left(\left(Y_{t}^{[L]}\right)_{t \leq t_n} \times \left(Y_{t}^{[R]}\right)_{t \leq t_n} = \left(y_{t}^{[L]}\right)_{t \leq t_n} \times \left(y_{t}^{[R]}\right)_{t \leq t_n}\right)\\
	&= q^{-m} (q)_{m}^2q^{m+\sum_{t \leq t_n} \left(y_{t}^{[L]}+y_{t}^{[R]}\right)}\sum_{\lambda} q^{|\lambda|} ,
\end{align*}
where the sum is taken over pairs of partitions $\lambda$ with parts at most $y_{t_n}^{[j]}$, respectively, for $j \in \{L,R\}$.  Continuing using Lemma \ref{L:largepartsproduct}, this is
\begin{align*}
	\frac{q^{\sum\limits_{t \leq t_n} \left(y_{t}^{[L]}+y_{t}^{[R]}\right)} (q)_{m}^2}{(q)_{y_{t_n}^{[L]}}(q)_{y_{t_n}^{[R]}}}
	&=q^{\sum\limits_{t \leq t_n} \left(y_{t}^{[L]}+y_{t}^{[R]}\right)} \prod_{y_{t_n}^{[L]} < t \leq m} \left(1-q^t\right) \prod_{y_{t_n}^{[R]} < t \leq m} \left(1-q^t\right) \\
	&=e^{-\sum\limits_{t \leq t_n} \left(w^{[L]}_t+w^{[R]}_t\right)} \left(\frac{1}{2B\sqrt{n}}\right)^{2t_n} \hspace{-.2cm}\prod_{y_{t_n}^{[L]} < t \leq m} \left(1-q^t\right) \prod_{y_{t_n}^{[R]} < t \leq m} \left(1-q^t\right)\\
	&\sim e^{-\sum\limits_{t \leq t_n} \left(w^{[L]}_{t_n}+w^{[R]}_{t_n}\right)} \left(\frac{1}{2B\sqrt{n}}\right)^{2t_n} e^{e^{-r}-\frac{e^{-w^{[L]}_{t_n}}}{2}-\frac{e^{-w^{[R]}_{t_n}}}{2}} \\
	&= \frac{1}{2^{2t_n}}   e^{e^{-r}-\sum\limits_{t \leq t_n} (w^{[L]}_{t}+w^{[R]}_{t})-\frac{e^{-w^{[L]}_{t_n}}}{2}-\frac{e^{-w^{[R]}_{t_n}}}{2}} \left(\frac{1}{B\sqrt{n}}\right)^{2t_n},
\end{align*}
 uniformly for $r, w^{[L]}_t,w^{[R]}_t \geq -\frac{\log(n)}{8}$.  Now let
$$
S:=\left\{{\bm w} : w^{[L]}_{t_n} \geq - \frac{\log (n)}{8} \ \text{and} \ w^{[R]}_{t_n} \geq - \frac{\log (n)}{8}  \right\}.
$$
Since $w^{[j]}_t \in \frac{1}{B\sqrt{n}}(\mathbb{Z}-\log(2B\sqrt{n}))$, recognizing Riemann sums gives
$$
\bm{\zeta}_{n,m}(U \cap S) \sim \bm{\nu}_{n,r}(U \cap S),
$$
for $U$ as in \eqref{E:openUdef} uniformly for $r$ in any $[r_1,r_2]$; in particular, we have $\bm{\zeta}_{n,m}(S) \sim \bm{\nu}_{n,r}(S)$.  But, letting $S^c$ denote the complement of $S$, we see that $\bm{\nu}_{n,r}(S^c) \to 0$ follows exactly as in \cite{F}, p. 724.  Now note that
$$
0=1-1=\bm{\zeta}_{n,m}(S^{c})-\bm{\nu}_{n,r}(S^c)+\bm{\zeta}_{n,m}(S)-\bm{\nu}_{n,r}(S) = \bm{\zeta}_{n,m}(S^c) +o(1).
$$  Thus, $\bm{\zeta}_{n,m}(S^c) \to 0$ also, and we have $\bm{\zeta}_{n,m}(U) \sim \bm{\nu}_{n,r}(U)$ uniformly for $r$ in any $[r_1,r_2]$, as desired.

Recalling that
$$
Y_{t}^{[j]}=\sup\left\{\ell : \sum_{k \geq \ell} X_k^{[j]} \geq t\right\},
$$
we see that $W_n \in S$ if and only if for $j \in \{L,R\}$
$$
\sup\left\{\ell : \sum_{k \geq \ell} X_k^{[j]} \geq t_n\right\} \geq B\sqrt{n}\log\left(2Bn^{\frac{3}{8}}\right),
$$
i.e., if and only if $Y_{t_n}^{[j]}$ depends only on $X_k^{[j]}$ for $k \geq B\sqrt{n}\log(2Bn^{\frac{3}{8}})$.  But now taking $a_n=B\sqrt{n}\log(2Bn^{\frac{3}{8}})$ and $b_n=n$, we have
$$
\sum_{a_n \leq k \leq b_n} \frac{k^2q^k}{\left(1-q^k\right)^2} \leq B^3n^{\frac 32} \int_{\log\left(2Bn^{\frac 38}\right)}^{\infty} \frac{u^2e^{-u}}{(1-e^{-u})^2}du = o\left(n^{\frac 32}\right),
$$
and Proposition \ref{P:W_ncondition} applies to $\bm{X}_{[a_n,n]}$.  This in turn yields $\bm{\xi}_{n,m}(U \cap S) \sim \bm{\zeta}_{n,m}(U \cap S)$ for $U$ as above, and then finally $\bm{\xi}_{n,m}(U) \sim \bm{\zeta}_{n,m}(U)$ as before.  Using the uniformity of these estimates for $r \in [r_1,r_2]$ finishes the proof.

The proof of Theorem \ref{T:*largesmallparts} is the same, except we use Lemma \ref{L:largepartsproduct} to estimate the product that arises in the calculation of $\mathbf{Q}_{q,m}^*(W_n={\bm w})$. \hspace{10cm} \qedsymbol

\subsection{Total small parts: proof of Theorem \ref{T:totalsmallparts}}

As in Subsection 5.1, Proposition \ref{P:W_ncondition} applies to $\bm{X}_{[1,k_n]}$ with $k_n=o(n^{\frac{1}{2}})$.  This clearly implies that for 
$$
W_n:= \left(\frac{\sum_{k \leq k_n} X_{k}^{[j]}-B\sqrt{n}k\log \left(k_n\right)}{B\sqrt{n}} \right)_{j \in \{L,R\}},
$$
we have $d_{\rm{TV}}(\mathbf{P}_{n,m}(W_n^{-1}), \mathbf{Q}_{q,m}(W_n^{-1})) \to 0.$  Now, as in Section \ref{S:smallparts}, we write
\begin{multline}\label{E:totalsmpartsmsummand}
\mathbf{P}_{n}\left(\left(W_n\right)_{j} \leq v_{j}, \ j \in \{L,R\}\right) \\ 
=\left(\sum_{r < r_1} + \sum_{r \in [r_1,r_2]} + \sum_{r > r_2} \right)\mathbf{P}_{n,m}\left(\left(W_n\right)_j \leq v_{j}, \ j \in \{L,R\}\right)\mathbf{P}_{n}\left({\rm PK}=m\right).
\end{multline}
We ignore the sums over the ranges $r<r_1$ and $r>r_2$ which  tend to 0, and in the range $r \in [r_1,r_2]$ we may replace $\mathbf{P}_{n,m}$ with $\mathbf{Q}_{q,m}$, so that (using independence) this is asymptotic to
\begin{equation}\label{E:totalsmpartsQreplace}
\sum_{r \in [r_1,r_2]} \mathbf{P}_{n,m}({\rm PK}=m) \prod_{j \in \{L,R\}}\mathbf{Q}_{q,m}\left(\left(W_n\right)_{j} \leq v_j\right).
\end{equation}

Following \S 8 of \cite{F}, we focus now on a particular term of $\mathbf{Q}_{q,m}$ and first restrict the range to $k \leq \ell_n$, where $\ell_n:=\lfloor k_n^{\frac 12} \rfloor = o(n^{\frac 14})$, so that we can use the calculation \eqref{E:smpartsQprobs}.  It is also simpler if we first do this without subtracting the $B\sqrt{n}\log \left(k_n\right)$ term.  Thus, we have
\begin{multline}
\mathbf{Q}_{q,m}\left(\frac{\sum_{k \leq \ell_n} X_{k}^{[j]}}{B\sqrt{n}} \leq v_j \right)  \\
=\sum_{w^{[j]}_{k_n} \in \frac{1}{B\sqrt{n}}\mathbb{N}_0 \cap [0,v_j]} \cdots \sum_{w^{[j]}_1 \in \frac{1}{B\sqrt{n}}\mathbb{N}_0 \cap \left[0,v_j-w^{[j]}_{\ell_n}-\dots - w^{[j]}_{2}\right]} \prod_{k=1}^{\ell_n} \mathbf{Q}_{q,m}\left(\frac{X_{k}^{[j]}}{B\sqrt{n}} = w^{[j]}_k\right). \label{E:totalsmpartsRiemannsum}
\end{multline}

Now, by \eqref{E:smpartsQprobs},
$$
\prod_{k=1}^{\ell_n} \mathbf{Q}_{q,m}\left(\frac{X_{k}^{[j]}}{B\sqrt{n}} = w^{[j]}_k\right)=\prod_{k=1}^{\ell_n} \mathbf{Q}_{q,m}\left(X_{k}^{[j]} = \frac{B\sqrt{n}}{k}kw^{[j]}_k\right) \sim \ell_n! \prod_{k=1}^{\ell_n} \left(e^{-kw^{[j]}_k} \frac{1}{B\sqrt{n}}\right).
$$
Plugging this into \eqref{E:totalsmpartsRiemannsum} and recognizing Riemann sums, we get, by Lemma \ref{L:stochYuleintegral}
\begin{align*}
\mathbf{Q}_{q,m}\left(\frac{\sum_{k \leq \ell_n} X_{k}^{[j]}}{B\sqrt{n}} \leq v_j\right) &\sim \ell_n!\int_0^{v_j} \cdots \int_0^{v_j-u^{[j]}_{\ell_n}-\dots -u^{[j]}_2} e^{-u_{1}^{[j]}-2u^{[j]}_2-\dots \ell_nu^{[j]}_{\ell_n}}du_{1}^{[j]} \cdots du^{[j]}_{\ell_n} \\
&=\left(1-e^{-v_j}\right)^{\ell_n}.
\end{align*}
 Since this is uniform in $v_j \in [0, \infty)$, we may replace $v_j \mapsto v_j + \log (\ell_n)$  (for fixed $v_j$), to get
$$
\mathbf{Q}_{q,m}\left(\frac{\sum_{k \leq \ell_n} X_{k}^{[j]}-B\sqrt{n}\log \left(\ell_n\right)}{B\sqrt{n}} \leq v_j\right) \sim e^{-e^{-v_j}}.
$$
But we want to show that the above holds with $\ell_n \mapsto k_n$.  But since $\ell_n = \lfloor k_n^{\frac 12} \rfloor$, this is equivalent to proving that
$$
\frac{\sum_{k_n^{\frac 12} < k \leq k_n} X_{k}^{[j]}-B\sqrt{n}\log \left(k_n^{\frac 12}\right)}{B\sqrt{n}}
$$
asymptotically has a point-mass distribution with mean 0.  This is accomplished by showing that its expectation and variance under $\mathbf{Q}_{q,m}$ are both $o(1)$.  This in turn follows from
\begin{align*}
&\mathbf{Var}_{q,m}\left(\frac{\sum_{k_n^{\frac 12} < k \leq k_n} X_{k}^{[j]}-B\sqrt{n}\log \left(k_n^{\frac 12}\right)}{B\sqrt{n}} \right)
=\mathbf{Var}_{q,m}\left(\frac{1}{B\sqrt{n}} \sum_{k_n^{\frac 12} < k \leq k_n} X_{k}^{[j]}\right)\\
&= \frac{1}{B^2n}\sum_{k_n^{\frac 12} < k \leq k_n} \frac{q^k}{\left(1-q^k\right)^2} 
\sim \frac{1}{B^2n}\sum_{k_n^{\frac 12} < k \leq k_n} \frac{1}{\left(1-e^{-\frac{k}{B\sqrt{n}}}\right)^2} 
\sim \sum_{k_n^{\frac 12}<k \leq k_n} \frac{1}{k^2} =o(1),
\\
&\mathbf{E}_{q,m}\left(\frac{1}{B\sqrt{n}}\sum_{k_n^{\frac 12} < k \leq k_n} X_k^{[j]}\right)
=\frac{1}{B\sqrt{n}}\sum_{k_n^{\frac 12} < k \leq k_n} \frac{q^k}{1-q^k}= \sum_{k_n^{\frac 12} < k \leq k_n} \frac{1}{k} +o(1) = \log \left(k_n^{\frac 12}\right) +o(1).
\end{align*}

Thus, we have
$$
\mathbf{Q}_{q,m}\left(\left(W_n\right)_j \leq v_j \right) = \mathbf{Q}_{q,m}\left(\frac{\sum_{k \leq k_n} X_{k}^{[j]}-B\sqrt{n}\log \left(k_n\right)}{B\sqrt{n}} \leq v_j\right) \sim e^{-e^{-v_j}}.
$$
Using this in \eqref{E:totalsmpartsQreplace} and then replacing $\mathbf{P}_{n,m}$ in \eqref{E:totalsmpartsmsummand} with $\mathbf{Q}_{q,m}$, the range over $r \in [r_1,r_2]$ is asymptotic to
$$
\sum_{r \in [r_1,r_2]} \mathbf{P}_{n}({\rm PK}=m) e^{-e^{-v_{L}}-e^{-v_R}} \sim \left(e^{-e^{-r_2}}-e^{-e^{-r_1}}\right)e^{-e^{-v_L}-e^{-v_R}},
$$
by Proposition \ref{P:denomlocalasymp}.  Taking $r_2 \to \infty$ and $r_1 \to - \infty$ completes the proof. \qed

\section{Moment Generating Functions}\label{S:moments}
In this section, we outline an approach towards refining Theorem \ref{T:PKdist} with the method of moments.  This discussion is independent of our conditioned Boltzmann model. Throughout this section, we assume that $q \in \mathbb{C}$ with $|q|<1$.  Recall that $Y_1(\lambda)$ denotes the largest part in the integer partition $\lambda$.  By standard combinatorial arguments, one finds the generating function for the $k$-th moment of the largest part for partitions to be
$${\rm MP}_{k}(q)=\sum_{n \geq 0} {\rm mp}_{k}(n)q^n:=\sum_{\lambda \in \mathcal{P}} \lambda_1(\lambda)^kq^{|\lambda|} = \sum_{m \geq 0} \frac{m^kq^m}{(q)_m}, \qquad \text{for $k \geq 0$.}$$
Analogous to the mean found in Theorem \ref{T:PKdist}, Theorem \ref{T:ErdosLehner} implies that
$$
{\rm mp}_{1}(n) = A\sqrt{n}\log(A)\sqrt{n} + A \gamma \sqrt{n} \left(1+o(1)\right),
$$
where $A=\frac{\sqrt{6}}{\pi}$, as before.  Ngo and Rhoades used the factorization (\cite{NR}, equation (1.8)),
$$
{\rm MP}_{1}(q)={\rm MP}_{0}(q)\sum_{n \geq 1} \frac{q^n}{1-q^n}=\frac{1}{(q)_{\infty}}\sum_{n \geq 1} \frac{q^n}{1-q^n},
$$
essentially a product of a modular form and a quantum modular form, to improve the error term to $\log(n)$.
\begin{theorem}[Theorem 1.5 in \cite{NR}]  We have
$${\rm mp}_{1}(n) = A\sqrt{n}\log\left(A\sqrt{n}\right) + A \gamma \sqrt{n}+O\left(\log (n)\right).$$
\end{theorem}
Furthermore, they found the recursions (\cite{NR}, remark on p. 10)
\begin{equation}\label{E:partmomentrec}
{\rm MP}_{k}(q)=\sum_{j=0}^{k-1} \binom{k-1}{j}{\rm MP}_{j}(q)S_{k-1-j}(q), \qquad \text{where} \qquad S_{k}(q):=\sum_{n \geq 1} \frac{n^kq^n}{1-q^n}.
\end{equation}
which express each ${\rm MP}_{k}$ recursively in terms of modular forms and quantum modular forms.  They stated that their methods could be used to prove asymptotic expansions for all moments ${\rm mp}_{k}(n)$.

Turning to unimodal sequences, recall that $\mathrm{PK}(\lambda)$ denotes the size of the peak.  Let
$$
{\rm MU}_{k}(q)=\sum_{n \geq 0} {\rm mu}_{k}(n)q^n:=\sum_{\lambda \in \mathcal{U}} \mathrm{PK}(\lambda)^kq^{|\lambda|}=\sum_{m \geq 0} \frac{m^kq^m}{(q)_m^2}, \qquad \text{for $k \in \N_0$.}
$$
In particular, the generating function for unimodal sequences satisfies (see \cite{St}, Proposition 2.5.1)
$$
U(q)={\rm MU}_{0}(q)=\frac{1}{(q)_{\infty}^2}\sum_{n \geq 0} (-1)^{n}q^{\frac{n(n+1)}{2}},
$$
which is a product of a modular form and a false theta function.  Recently, Nazaroglu and the second author discovered how to fit this false theta function into a modular framework \cite{BN}.  Thus, it would be interesting if we could relate the higher moments ${\rm MU}_{k}$ to ${\rm MU}_{0}$ in the way that Ngo and Rhoades did for partitions.  We leave this as an open problem, but we prove a recurrence that is somewhat analogous to \eqref{E:partmomentrec}.

Recall that the complete Bell polynomials $\mathbb{B}_k=\mathbb{B}_k(a_1, \dots, a_k)$ are defined by
$$
\sum_{k \geq 0} \frac{\mathbb{B}_k}{k!}u^k:=\exp\left(\sum_{k \geq 1} \frac{a_k}{k!}u^k\right).
$$

\begin{theorem}\label{thm:5.1}

For $k \geq -1$ and $n\geq 1$, define
$$
S_{k,n}(q):=\sum_{m \geq 1} \frac{m^{k}q^{nm}}{1-q^m}, 
$$
and let $\mathbb{B}_{k,n}:= \mathbb{B}_{k}\left(S_{0,n}, \dots, S_{k-1,n}\right)$ be the complete Bell polynomials in $S_{j,n}$. Then
\begin{equation}\label{E:U_kgenfn}
{\rm MU}_{k}(q) = \frac{1}{(q)_{\infty}^2}\sum_{n \geq 0} (-1)^{n}q^{\frac{n(n+1)}{2}}\mathbb{B}_{k,n+1}.
\end{equation}
\end{theorem}

\begin{remark}  Note that we have the recurrence (see \cite{Co}, \S 3.3.)
$$
\mathbb{B}_k:=\begin{cases} 1 & \text{if $k=0$,} \\ \sum_{j=0}^{k-1} \left(\begin{smallmatrix}k-1\\j\end{smallmatrix}\right) \mathbb{B}_{k-j}a_{j+1} & \text{if $k \geq 1$.}\end{cases}
$$  But the reason that we cannot directly obtain a recurrence for the ${\rm MU}_k$s like in \eqref{E:partmomentrec} is that the $\mathbb{B}$s and $S$s depend on $n$.
\end{remark}

\begin{proof}[Proof of Theorem \ref{thm:5.1}]
Define
\begin{equation*}\label{E:U_k2vargenfn}
U(\zeta;q):=\sum_{\lambda \in \mathcal{U}}\zeta^{\mathrm{PK}(\lambda)}q^{|\lambda|}=\sum_{m \geq 0} \frac{\zeta^mq^m}{(q)_m^2}, 
\end{equation*} and let $\delta_{\zeta}:= \zeta\frac{\partial}{\partial\zeta}$.  Then 
\begin{equation}\label{E:M_kderiv}
{\rm MU}_{k}(q)= \left[\delta_{\zeta}^k \left(U(\zeta;q) \right) \right]_{\zeta=1}
\end{equation}
Using straightforward manipulation with Euler's two series expansions (\cite{A}, Corollary 2.2)
\begin{equation*}\label{E:Eulerdistinct}
\frac{1}{(\zeta;q)_{\infty}} = \sum_{n \geq 0} \frac{\zeta^n}{(q)_n}, \qquad (-\zeta;q)_{\infty}=\sum_{n \geq 0} \frac{\zeta^nq^{\frac{n(n-1)}{2}}}{(q)_n},
\end{equation*}
we have
\begin{align*}
U(\zeta;q)&=1 + \sum_{m \geq 1} \frac{\zeta^mq^m}{(q)_m^2} 
= 1 + \frac{1}{(q)_{\infty}}\sum_{m \geq 1} \frac{\zeta^{m}q^m\left(q^{m+1}\right)_{\infty}}{(q)_{m}} \\ 
&= 1 + \frac{1}{(q)_{\infty}}\sum_{\substack{m \geq 1 \\ n \geq 0}}\frac{(-1)^{n}\zeta^{m}q^{m+n(m+1) + \frac{n(n-1)}{2}}}{(q)_{m}(q)_{n}}  
= 1 + \frac{1}{(q)_{\infty}}\sum_{\substack{m \geq 1 \\ n \geq 0}}\frac{(-1)^{n}q^{\frac{n(n+1)}{2}}}{(q)_{n}}\cdot \frac{\zeta^mq^{(n+1) m}}{(q)_m} \\
&=1+\frac{1}{(q)_{\infty}}\sum_{n \geq 0} \frac{(-1)^{n}q^{\frac{n(n+1)}{2}}}{(q)_{n}}\left(\frac{1}{(\zeta q^{n+1})_{\infty}} - 1 \right).
\end{align*}
By the chain rule,
$$
\left[\delta_{\zeta}^k \left(\frac{1}{(\zeta q^{n+1})_{\infty}}\right) \right]_{\zeta=1} = \left[\frac{\partial^k}{\partial u^k}\frac{1}{(e^u q^{n+1})_{\infty}}\right]_{u=0},
$$
so it suffices to find the Taylor expansion of $(e^u q^{n+1})^{-1}_{\infty}$ about $u=0$.  We first find the Taylor expansion of its logarithm and then exponentiate, which introduces the complete Bell polynomials.  We have
$$
\Log \lp \frac{1}{(e^uq^{n+1})_{\infty}} \rp= -\sum_{\ell \geq n+1}  \Log\left(1-e^uq^\ell\right) = \sum_{\substack{\ell \geq n+1 \\ m \geq 1}} \frac{e^{mu}q^{m\ell}}{m} = \sum_{m \geq 1} \frac{e^{mu}q^{(n+1)m}}{m(1-q^m)} = \sum_{k \geq 0} \frac{u^k}{k!}S_{k-1,n+1},
$$
where clearly $S_{-1,n+1}=\Log (\frac{1}{(q^{n+1})_{\infty}})$.  Hence,
$$
\frac{1}{(e^uq^{n+1})_{\infty}}= \frac{1}{(q^{n+1})_{\infty}}\exp\left(\sum_{k \geq 1} \frac{u^k}{k!}S_{k-1,n+1}\right) = \frac{1}{(q^{n+1})_{\infty}}\sum_{k \geq 0} \frac{u^k}{k!} \mathbb{B}_{k,n+1}.
$$
Thus,
$$
\left[\delta_{\zeta}^k \left(\frac{1}{\left(\zeta q^{n+1}\right)_{\infty}}\right) \right]_{\zeta=1}=\frac{\mathbb{B}_{k,n+1}}{(q^{n+1})_{\infty}},
$$
and we obtain \eqref{E:U_kgenfn} by plugging into \eqref{E:M_kderiv}. 
\end{proof}
 
 We leave further exploration of the $q$-series ${\rm MK}_k(q)$ as an open problem.

\end{document}